\documentclass[11pt]{article}

\usepackage[T1]{fontenc}
\usepackage[utf8]{inputenc}
\IfFileExists{lmodern.sty}{\usepackage{lmodern}}{}
\usepackage[margin=1in]{geometry}
\usepackage{amsmath,amssymb,amsthm}
\usepackage{mathtools}
\usepackage{booktabs}
\usepackage{array}
\usepackage{longtable}
\usepackage{enumitem}
\usepackage[numbers,sort&compress]{natbib}
\usepackage[hidelinks]{hyperref}
\newcommand{\doi}[1]{\href{https://doi.org/\detokenize{#1}}{doi:\,\detokenize{#1}}}

\newtheorem{theorem}{Theorem}
\newtheorem*{theoremTwoPrime}{Theorem~2$'$}
\newtheorem*{theoremThreePrime}{Theorem~3$'$ (Rigidity)}
\newtheorem{lemma}{Lemma}[section]
\newtheorem{proposition}[lemma]{Proposition}
\newtheorem{corollary}[lemma]{Corollary}
\theoremstyle{definition}
\newtheorem{definition}[lemma]{Definition}
\newtheorem{example}[lemma]{Example}
\newtheorem{question}{Question}
\theoremstyle{remark}
\newtheorem{remark}[lemma]{Remark}

\newcommand{\R}{\mathbb{R}}
\newcommand{\Z}{\mathbb{Z}}
\newcommand{\Q}{\mathbb{Q}}
\newcommand{\N}{\mathbb{N}}
\newcommand{\C}{\mathbb{C}}
\newcommand{\T}{\mathbb{T}}
\newcommand{\Fib}{\mathrm{Fib}}
\newcommand{\bphi}{\bar{\varphi}}
\newcommand{\nint}[1]{\lVert #1 \rVert}
\newcommand{\Dl}{\Delta}
\newcommand{\Dset}{\mathcal{D}}
\newcommand{\Kp}{K_{\mathrm{P}}}
\DeclareMathOperator{\sgn}{sgn}

\title{The Golden Quadruple:\\[0.3em]
\large Four Reals Whose Differences Form a Geometric Progression,\\
Rigidity of the Golden Ratio, and a Bohr Almost Periodic Cosine Sum}
\author{Micha{\l} Przybylski, M.Sc.\ Eng.\thanks{ORCID: \href{https://orcid.org/0009-0005-5170-917X}{https://orcid.org/0009-0005-5170-917X}.}}
\date{Preprint, version 3, November 14, 2026}

\begin{document}
\maketitle

\begin{abstract}
Let $\varphi=(1+\sqrt5)/2$. The set $\{0,\varphi,\varphi^2,\varphi^3\}$ is classical (Euclid, \emph{Elements}~XIII.8):
by $\varphi^2=\varphi+1$ its six pairwise differences take the values $1,\varphi,\varphi^2,\varphi^3$, with multiplicities $1,2,2,1$. We
write it as the \emph{golden quadruple} $S=(F,\,F-a,\,F-b,\,F-a-b)$, $a=\varphi u$, $b=\varphi^2u$, $F>a+b$, and study its
extremal and rigidity properties. No five distinct reals have all their positive differences in $\{u\varphi^n:n\in\Z\}$; this follows from
a result of Lenstra (1977), and we give a self-contained proof. We give an elementary proof of uniqueness: every four-element set with
this property is a translate of $u\varphi^m\{0,\varphi,\varphi^2,\varphi^3\}$ for some $m\in\Z$, and if all its differences lie in
$\{u,\varphi u,\varphi^2u,\varphi^3u\}$ it is a translate of $S$. Among all $x>1$, only $x=\varphi$ admits a four-point set of reals whose set of distinct
differences is $\{c,cx,cx^2,cx^3\}$. For the analogous quadruple with base $x\neq\varphi$ the failure is measured by two closure residues,
which along Fibonacci quotients equal $\pm u/\Fib_n^2$ and $\pm \Fib_{n+1}u/\Fib_n^3$. The cosine sum
$C_S(t)=\sum\cos(2\pi d t)$ over the six differences equals $8\cos^2(\pi a t)\cos^2(\pi b t)-2$; it is not periodic
but is Bohr almost periodic, with explicit almost periods $\tau_n=\Fib_n/u$, and we determine the two-term asymptotics
$\sup_t|C_S(t+\tau_n)-C_S(t)|=6\sqrt3\,\pi\varphi^{-n}-2\sqrt3\,\pi^3\varphi^{-3n}+O(\varphi^{-5n})$, with optimal constants; the
inclusion length of the $\varepsilon$-almost periods is of order $1/\varepsilon$. For a scale $K>0$ the golden ladder
$\nu_n=K\varphi^{-n}$ realises the construction with differences $\nu_{200},\dots,\nu_{203}$; we evaluate it for
$K=\sqrt{c^5/(Gh)}$. All proofs are elementary and self-contained.
\end{abstract}



\tableofcontents

\section{Introduction}\label{sec:intro}

For a finite set $A\subset\R$, the multiset of its positive differences records a large part of
its additive structure. Two classical extremes are well understood. At one end are
\emph{Sidon sets} and \emph{Golomb rulers}, in which all positive differences are distinct
\cite{Sidon1932,ErdosTuran1941,OBryant2004,Golomb1972,Drakakis2009}. At the other end are
arithmetic progressions, whose difference set $A-A$ is as small as possible,
$|A-A|=2|A|-1$ \cite{Nathanson1996,TaoVu2006}. Between them lie the $(v,k,\lambda)$-difference
sets of design theory, in which every non-zero group element is represented the same number
of times \cite{Singer1938,Baumert1971,BJL1999}.

In this paper we impose a regularity of a different kind: we prescribe the \emph{shape} of the
set of distinct differences and ask that it be a geometric progression. The simplest
non-trivial case already singles out the golden ratio $\varphi=(1+\sqrt5)/2$. Fix $u>0$ and put
\[
  a=\varphi u,\qquad b=\varphi^2u,\qquad F>a+b,\qquad
  S=(F,\;F-a,\;F-b,\;F-a-b).
\]
We call $S$ the \emph{golden quadruple}. Because $\varphi^2=\varphi+1$, the identities
$b-a=u$ and $a+b=\varphi^3u$ hold simultaneously, and this is the whole mechanism behind the
first result.

The configuration itself is classical. As a set, $S$ is a translate of $u\{0,\varphi,\varphi^2,\varphi^3\}$
(Remark~\ref{rem:structure}(a)), and $\{0,\varphi,\varphi^2,\varphi^3\}$ is, up to a dilation, the set formed by
the two endpoints and the two crossing points on a side of the regular pentagram: by Euclid's
\emph{Elements}~XIII.8, two diagonals of a regular pentagon that subtend consecutive angles cut one another in
extreme and mean ratio, and their greater segments are equal to the side of the pentagon \cite{Heath1926}. The
details are recalled in Remark~\ref{rem:euclid}. Up to dilation it is also the sequence $0,1,\varphi,\varphi^2$ that Lenstra obtains in \cite[(3.9)]{Lenstra1977} from \cite[(2.4)(b)]{Lenstra1977}. Neither the set nor its difference multiset
(Theorem~\ref{thm:T1}) is new, and neither is the bound of four points
(Remark~\ref{rem:lenstra}); what this paper adds are the rigidity theorem, an explicit geometric formulation of uniqueness with an elementary proof, and the analysis of the
cosine sum.

\paragraph{Main results.} Write $\Dset=\{u,\varphi u,\varphi^2u,\varphi^3u\}$.
\begin{enumerate}[label=\textbf{(\arabic*)},leftmargin=2.2em]
\item \textbf{Theorem~\ref{thm:T1}} (classical). The six pairwise differences of $S$ are exactly
  $u,\varphi u,\varphi^2u,\varphi^3u$, with multiplicities $1,2,2,1$.
\item \textbf{Theorem~\ref{thm:T2} and Theorem~2$'$.} Every four-element set with all its positive differences
  in $\Dset$, or even in $\{c\varphi^n:n\in\Z\}$ for a fixed $c>0$, is a translate of a golden quadruple (an explicit geometric formulation with our own elementary proof; it also follows from Nagell's results via normalisation, Remark~\ref{rem:nagell}).
  No set of five or more distinct reals has this property; this bound follows from a result of
  Lenstra~\cite{Lenstra1977} (Remark~\ref{rem:lenstra}), and we include a self-contained elementary proof.
\item \textbf{Theorem~\ref{thm:T3}.} For a base $x>1$, $x\neq\varphi$, the analogous quadruple
  $(F,F-xu,F-x^2u,F-xu-x^2u)$ does not have the four consecutive powers $u,xu,x^2u,x^3u$ as its
  differences (nor any four-term geometric progression); the defects are the closure residues
  $\rho_-(x)=(x^2-x-1)u$ and $\rho_+(x)=x(x^2-x-1)u$, computed explicitly along Fibonacci quotients.
\item \textbf{Theorem~3$'$.} If $x>1$ and some four-element $T\subset\R$ satisfies
  $\Dl(T)=\{c,cx,cx^2,cx^3\}$, then $x=\varphi$ and $T$ is a translate of $c\{0,\varphi,\varphi^2,\varphi^3\}$.
\item \textbf{Theorem~\ref{thm:T4}.} The cosine sum
  \[
    C_S(t)=\cos(2\pi ut)+2\cos(2\pi\varphi ut)+2\cos(2\pi\varphi^2ut)+\cos(2\pi\varphi^3ut)
  \]
  is not periodic but is almost periodic in the sense of Bohr \cite{Bohr1925}. Moreover,
  \[
    C_S(t)=8\cos^2(\pi at)\cos^2(\pi bt)-2,
  \]
  and the numbers $\Fib_n/u$ are explicit almost periods. Proposition~\ref{prop:sharp} gives the sharp constant:
  $\sup_t|C_S(t+\Fib_n/u)-C_S(t)|=6\sqrt3\,\pi\varphi^{-n}-2\sqrt3\,\pi^3\varphi^{-3n}+O(\varphi^{-5n})$, and
  Proposition~\ref{prop:Leps} shows $L(\varepsilon)\asymp1/\varepsilon$ for the almost periods.
\end{enumerate}

\paragraph{Scale.} All statements are invariant under dilations, so the only role of a scale is
to fix a normalisation. For an arbitrary $K>0$ we use the \emph{golden ladder}
$\nu_n=K\varphi^{-n}$ ($n\in\Z$), which satisfies $\nu_n=\nu_{n+1}+\nu_{n+2}$. Choosing $u=\nu_{203}$
gives $a=\nu_{202}$, $b=\nu_{201}$, $a+b=\nu_{200}$, and the difference set of $S$ becomes
$\{\nu_{200},\nu_{201},\nu_{202},\nu_{203}\}$ (Corollary~\ref{cor:ladder}). As a single numerical
example we take $K=\sqrt{c^5/(Gh)}$ (Example~\ref{ex:planck}); it is used only as a positive real
number.

\paragraph{What is classical and what is new.} Everything here is elementary. The set
$\{0,\varphi,\varphi^2,\varphi^3\}$ goes back to Euclid (\emph{Elements}~XIII.8), and Theorem~\ref{thm:T1} is an
immediate consequence of $\varphi^2=\varphi+1$; we include its proof only to fix notation. The contributions of this paper are as follows; only the last two items are claimed as new
(Section~\ref{sec:literature}).
\begin{itemize}
\item Theorems~\ref{thm:T2}(ii) and~2$'$(ii) (uniqueness; not claimed as new): an explicit geometric formulation with our
  own elementary proof; they can also be derived from Nagell's list of exceptional units via normalisation
  \cite{LeutbecherMartinet1982,Nagell1964,Nagell1969a,Nagell1969b}, see Remark~\ref{rem:nagell}. The four-point
  configuration is unique up to translation (and, in Theorem~2$'$(ii), dilation by a power of~$\varphi$);
\item \emph{New:} Theorem~3$'$ (rigidity): among all ratios $x>1$, only $x=\varphi$ admits a four-point set whose set of
  distinct differences is $\{c,cx,cx^2,cx^3\}$; an explicit geometric rigidity statement with a self-contained elementary
  proof (it can also be derived from \cite[Lemma~3]{CampolinaMailybaev2021} by a short argument, Remark~\ref{rem:CM});
\item \emph{New:} the quantitative analysis of the cosine sum $C_S$: the exact norm of the shift by the almost periods
  $\Fib_n/u$ and its sharp asymptotics (Proposition~\ref{prop:sharp}), and the order $L(\varepsilon)\asymp1/\varepsilon$
  with explicit constants (Proposition~\ref{prop:Leps}). They rest on the closed form (Proposition~\ref{prop:closed}) and
  the almost periods $\Fib_n/u$ (Proposition~\ref{prop:fibAP}), which are elementary; the explicit value $C_S(\Fib_n/u)$ in Proposition~\ref{prop:fibAP}(iii) is ours. Theorem~\ref{thm:T4}
  (non-periodicity and Bohr almost periodicity) is a direct consequence of classical facts and is included for completeness.
\end{itemize}
The bound of four points (Theorems~\ref{thm:T2}(i) and~2$'$(iii)), i.e.\ $N(\varphi)=4$ in the notation of
Definition~\ref{def:N}, is \emph{not} new: it follows from Lenstra's work on Euclidean number fields
\cite[(2.2), (3.9)]{Lenstra1977}, see Remark~\ref{rem:lenstra}. Version~1 of this preprint listed it among the new
results; we correct this here. Our proof is elementary and independent of~\cite{Lenstra1977}. Likewise
Lemma~\ref{lem:unit} and Theorem~2$'$(i) follow from Nagell's determination of the exceptional units of
$\Q(\sqrt5)$ \cite{Nagell1964,Nagell1969a,Nagell1969b}.
Theorem~\ref{thm:T3} and Proposition~\ref{prop:fibres} are elementary companions of Theorem~3$'$ for quadruples
of the special form $S_x$. Section~\ref{sec:literature} discusses the relation to the literature in detail; it also
shows that some of these phenomena are genuinely specific to $\varphi$: for the plastic number there are
five-point sets all of whose differences are powers of the base.

\paragraph{Organisation.} Section~\ref{sec:prelim} fixes notation. Sections~\ref{sec:T1}--\ref{sec:T4}
contain the main theorems and their strengthenings. Section~\ref{sec:numerics} gives numerical
examples and tables, Section~\ref{sec:literature} discusses the literature and open questions, and
Section~\ref{sec:formal} describes which of the results have been formally verified in Lean~4.

\section{Preliminaries}\label{sec:prelim}

Throughout, $\N=\{0,1,2,\dots\}$, $\nint{\theta}=\min_{m\in\Z}|\theta-m|$ denotes the distance from
$\theta\in\R$ to the nearest integer, and $e(\theta)=e^{2\pi i\theta}$.
The Fibonacci numbers are $\Fib_0=0$, $\Fib_1=1$, $\Fib_{n+1}=\Fib_n+\Fib_{n-1}$.

\subsection{Differences}

\begin{definition}\label{def:diff}
Let $A\subset\R$ be finite. The \emph{representation function} of $A$ is
\[
  r_A(d)=\#\{(x,y)\in A\times A:\ x-y=d\}\qquad(d>0),
\]
the \emph{positive difference set} is $\Dl(A)=\{d>0:\ r_A(d)\ge1\}$, and the
\emph{difference multiset} of $A$ is the multiset in which each $d\in\Dl(A)$ occurs $r_A(d)$ times;
it has $\binom{|A|}{2}$ elements. We also use the full difference set $A-A=\{x-y:x,y\in A\}$.
For an ordered tuple we use the same notions for the set of its entries.
\end{definition}

\begin{definition}[classical notions, used in Section~\ref{sec:literature}]\label{def:classical}
A finite $A\subset\R$ is a \emph{Sidon set} if $r_A(d)\le1$ for all $d>0$. A \emph{Golomb ruler} is a
Sidon subset of $\Z$; it is \emph{perfect} if $\Dl(A)=\{1,2,\dots,\max A-\min A\}$. If $G$ is a finite
abelian group of order $v$, a $k$-subset $B\subseteq G$ is a $(v,k,\lambda)$-\emph{difference set} if
every $g\in G\setminus\{0\}$ has exactly $\lambda$ representations $g=x-y$ with $x,y\in B$.
\end{definition}

\begin{lemma}\label{lem:nint}
For $\theta,\theta'\in\R$ and $m\in\Z$:
$\nint{\theta+m}=\nint{\theta}=\nint{-\theta}$, $\nint{\theta+\theta'}\le\nint{\theta}+\nint{\theta'}$, and
$|e(\theta)-1|=2|\sin\pi\theta|\le2\pi\nint{\theta}$.
\end{lemma}

\begin{proof}
The first two identities are clear. If $m,m'\in\Z$ are nearest to $\theta,\theta'$, then
$\nint{\theta+\theta'}\le|\theta+\theta'-m-m'|\le\nint{\theta}+\nint{\theta'}$.
Next, $|e(\theta)-1|=|e^{\pi i\theta}-e^{-\pi i\theta}|=2|\sin\pi\theta|$, and for $m$ nearest to
$\theta$ we have $|\sin\pi\theta|=|\sin\pi(\theta-m)|\le\pi|\theta-m|=\pi\nint{\theta}$.
\end{proof}

\subsection{The golden ratio}

Let $\varphi=(1+\sqrt5)/2$ and $\bphi=(1-\sqrt5)/2$.

\begin{lemma}\label{lem:phi}
\begin{enumerate}[label=\textup{(\roman*)}]
\item $\varphi^2=\varphi+1$, $\varphi+\bphi=1$, $\varphi\bphi=-1$, hence $\bphi=-\varphi^{-1}$.
\item $\varphi^{n+2}=\varphi^{n+1}+\varphi^n$ for every $n\in\Z$. In particular
  $\varphi^3=2\varphi+1$, $\varphi^4=3\varphi+2$, $\varphi-1=\varphi^{-1}$ and $2-\varphi=\varphi^{-2}$.
\item $3/2<\varphi<2$; consequently $\varphi^2<3$ and $\varphi^3>4$.
\item $\varphi$ is irrational; equivalently, $1$ and $\varphi$ are linearly independent over $\Q$.
\item For $x>1$, $x^2-x-1=(x-\varphi)(x-\bphi)$ with $x-\bphi>0$; hence $x^2-x-1=0$ if and only if
  $x=\varphi$, and $\sgn(x^2-x-1)=\sgn(x-\varphi)$.
\item The map $n\mapsto\varphi^n$ is strictly increasing on $\Z$, in particular injective.
\end{enumerate}
\end{lemma}

\begin{proof}
(i) $\varphi^2=(6+2\sqrt5)/4=(3+\sqrt5)/2=\varphi+1$; the other two identities are immediate, and
$\bphi=-1/\varphi$ follows from $\varphi\bphi=-1$.
(ii) Multiply $\varphi^2=\varphi+1$ by $\varphi^n$. With $n=1,2$ we get $\varphi^3=\varphi^2+\varphi=2\varphi+1$
and $\varphi^4=\varphi^3+\varphi^2=3\varphi+2$; with $n=-1$, $\varphi=1+\varphi^{-1}$; with $n=-2$,
$1=\varphi^{-1}+\varphi^{-2}$, so $\varphi^{-2}=1-(\varphi-1)=2-\varphi$.
(iii) From $2<\sqrt5<3$. Then $\varphi^2=\varphi+1<3$ and $\varphi^3=2\varphi+1>4$.
(iv) If $\varphi=p/q\in\Q$ then $\sqrt5=2p/q-1\in\Q$. But if $\sqrt5=m/n$ in lowest terms, then
$m^2=5n^2$, so $5\mid m$, hence $25\mid 5n^2$, so $5\mid n^2$ and $5\mid n$, a contradiction. If $r+s\varphi=0$ with
$r,s\in\Q$ not both zero, then $s\ne0$ and $\varphi=-r/s\in\Q$.
(v) The roots of $x^2-x-1$ are $\varphi$ and $\bphi<0<1<x$.
(vi) $\varphi>1$.
\end{proof}

\subsection{The construction}

\begin{definition}[golden quadruple]\label{def:S}
Fix $u>0$, put $a=\varphi u$, $b=\varphi^2u$, and fix $F>a+b$. The \emph{golden quadruple} is
\[
  S=S(F,u)=(s_1,s_2,s_3,s_4)=(F,\;F-a,\;F-b,\;F-a-b).
\]
We write $\Dset=\Dset(u)=\{u,\varphi u,\varphi^2u,\varphi^3u\}$.
Since $0<a<b<a+b<F$, we have $s_1>s_2>s_3>s_4>0$; the hypothesis $F>a+b$ ensures $S\subset(0,\infty)$.
\end{definition}
The name should not be confused with the \emph{golden quadruples} of Wei\ss{} \cite[Theorem~3]{Weiss2002}: four
collinear points characterised projectively by their cross ratios. That notion is projectively invariant and covers a
wider class than the metric condition of Definition~\ref{def:S}.

\begin{definition}[golden ladder]\label{def:ladder}
For $K>0$ and $n\in\Z$ put $\nu_n=\nu_n(K)=K\varphi^{-n}$.
\end{definition}

\begin{lemma}\label{lem:ladder}
For every $n\in\Z$: $\nu_n=\nu_{n+1}+\nu_{n+2}$ and $\nu_n=\varphi\,\nu_{n+1}$. If $u=\nu_{203}$, then
$a=\nu_{202}$, $b=\nu_{201}$, $a+b=\nu_{200}$ and $b-a=\nu_{203}$.
\end{lemma}

\begin{proof}
$\nu_n-\nu_{n+1}-\nu_{n+2}=K\varphi^{-n-2}(\varphi^2-\varphi-1)=0$, and $\nu_n/\nu_{n+1}=\varphi$. Hence
$a=\varphi\nu_{203}=\nu_{202}$ and $b=\varphi^2\nu_{203}=\nu_{201}$; the recursion with $n=200$ gives
$a+b=\nu_{202}+\nu_{201}=\nu_{200}$, and with $n=201$ gives $b-a=\nu_{201}-\nu_{202}=\nu_{203}$.
\end{proof}

\begin{example}[the Planck scale]\label{ex:planck}
Let $\Kp=\sqrt{c^5/(Gh)}$, where (in this example only) $c=299\,792\,458$, $G=6.674\,30\times10^{-11}$ and
$h=6.626\,070\,15\times10^{-34}$ are the CODATA~2018 numerical values (SI units) of the speed of light,
the Newtonian constant of gravitation and the Planck constant \cite{CODATA2018}; $c$ and $h$ are exact,
and $G$ has relative standard uncertainty $2.2\times10^{-5}$. Then
\[
  \Kp=7.399\,82(8)\times10^{42},
\]
the relative uncertainty $1.1\times10^{-5}$ being inherited from $G$. Since $\varphi^{203}\approx2.657\,61\times10^{42}$,
the index $203$ brings the ladder to order one: $\nu_{203}(\Kp)\approx2.784\,39$.
We keep $h$ as in the definition; with $\hbar=h/(2\pi)$ in its place one would obtain
$\sqrt{2\pi}\,\Kp\approx1.854\,86\times10^{43}$.
No interpretation of $\Kp$ beyond its numerical value is intended.
\end{example}

\subsection{Almost periodic functions and cosine sums}

\begin{definition}[Bohr \cite{Bohr1925,Bohr1947,Besicovitch1932,Corduneanu1989}]\label{def:ap}
A set $E\subseteq\R$ is \emph{relatively dense} if there is $L>0$ such that every closed interval of
length $L$ meets $E$. For a continuous $f:\R\to\C$ and $\varepsilon>0$, the set of
\emph{$\varepsilon$-almost periods} is
\[
  E(f;\varepsilon)=\Bigl\{\tau\in\R:\ \sup_{t\in\R}|f(t+\tau)-f(t)|\le\varepsilon\Bigr\}.
\]
The function $f$ is \emph{almost periodic} (in the sense of Bohr) if $E(f;\varepsilon)$ is relatively
dense for every $\varepsilon>0$. It is \emph{periodic} if $f(t+T)=f(t)$ for all $t$ and some $T>0$.
\end{definition}

A continuous periodic function with period $T$ is almost periodic, because $E(f;\varepsilon)\supseteq T\Z$
and every closed interval of length $T$ meets $T\Z$.

\begin{definition}[cosine sum]\label{def:cosine}
For a finite $A\subset\R$ put
\[
  C_A(t)=\sum_{\substack{x,y\in A\\ x>y}}\cos\bigl(2\pi(x-y)t\bigr)=\sum_{d\in\Dl(A)}r_A(d)\cos(2\pi dt),
  \qquad \widehat A(t)=\sum_{x\in A}e(xt).
\]
\end{definition}

\section{The difference multiset of the golden quadruple}\label{sec:T1}

\begin{lemma}[two-dimensional cubes]\label{lem:cube}
Let $0<\alpha<\beta$ and $Q=\{0,\alpha,\beta,\alpha+\beta\}=\{0,\alpha\}+\{0,\beta\}$. The six positive
differences of $Q$ are
\[
  \alpha-0=\alpha,\quad \beta-0=\beta,\quad (\alpha+\beta)-0=\alpha+\beta,\quad
  \beta-\alpha,\quad (\alpha+\beta)-\alpha=\beta,\quad (\alpha+\beta)-\beta=\alpha,
\]
so the difference multiset of $Q$ is $\{\beta-\alpha,\alpha,\alpha,\beta,\beta,\alpha+\beta\}$. For every
$F\in\R$ the reflected set $F-Q$ and every translate $t+Q$ have the same difference multiset.
Moreover $Q$ is symmetric: $(\alpha+\beta)-Q=Q$.
\end{lemma}

\begin{proof}
The list is a direct enumeration of the $\binom42=6$ pairs. For $x>y$ in $Q$ one has
$(F-y)-(F-x)=x-y$ and $(t+x)-(t+y)=x-y$, so reflection and translation preserve the difference
multiset. The symmetry is $\{\alpha+\beta,\beta,\alpha,0\}=Q$.
\end{proof}

\begin{theorem}\label{thm:T1}
The six pairwise differences of the golden quadruple are
\[
\begin{aligned}
  s_1-s_2&=\varphi u, & s_1-s_3&=\varphi^2u, & s_1-s_4&=\varphi^3u,\\
  s_2-s_3&=u,         & s_2-s_4&=\varphi^2u, & s_3-s_4&=\varphi u.
\end{aligned}
\]
Consequently $\Dl(S)=\Dset=\{u,\varphi u,\varphi^2u,\varphi^3u\}$, these four numbers are distinct and form
a geometric progression with ratio $\varphi$, and
\[
  r_S(u)=1,\qquad r_S(\varphi u)=2,\qquad r_S(\varphi^2u)=2,\qquad r_S(\varphi^3u)=1.
\]
\end{theorem}

\begin{proof}
Directly: $s_1-s_2=a=\varphi u$, $s_1-s_3=b=\varphi^2u$, $s_2-s_4=b$, $s_3-s_4=a$. Further
$s_1-s_4=a+b=(\varphi+\varphi^2)u=\varphi(1+\varphi)u=\varphi\cdot\varphi^2u=\varphi^3u$ and
$s_2-s_3=b-a=(\varphi^2-\varphi)u=u$, both by $\varphi^2=\varphi+1$ (Lemma~\ref{lem:phi}(i)). Since
$\varphi>1$, $u<\varphi u<\varphi^2u<\varphi^3u$, so the four values are distinct; counting occurrences
gives the multiplicities. (Equivalently: $S=F-Q$ with $Q=\{0,a\}+\{0,b\}$, and Lemma~\ref{lem:cube}
applies with $\beta-\alpha=u$ and $\alpha+\beta=\varphi^3u$.)
\end{proof}

\begin{remark}[classical origin of the configuration]\label{rem:euclid}
The set $u\{0,\varphi,\varphi^2,\varphi^3\}$, of which $S$ is a translate (Remark~\ref{rem:structure}(a)), is
classical. Let $ABCDE$ be a regular pentagon whose diagonals have length $\varphi^3u$, and let the diagonal $AC$
meet $BE$ at $H$ and $BD$ at $H'$. The diagonals $AC$ and $BE$ subtend the consecutive angles at $B$ and $A$, so by
Euclid's \emph{Elements}~XIII.8 \cite[Book~XIII, Prop.~8]{Heath1926} the point $H$ cuts $AC$ in extreme and mean
ratio, and the greater segment $CH$ is equal to the side of the pentagon, which is $\varphi^3u/\varphi=\varphi^2u$.
Hence $AH=\varphi^3u-\varphi^2u=\varphi u$. The reflection in the perpendicular bisector of $AC$ maps the pentagon
onto itself, exchanges $A$ with $C$ and $BE$ with $BD$, and therefore $H$ with $H'$; so $CH'=AH=\varphi u$ and
$AH'=\varphi^2u$. Measured from $A$, the four points $A,H,H',C$ of the segment $AC$, which is a side of the
pentagram $ACEBD$, lie at distances $0,\ \varphi u,\ \varphi^2u,\ \varphi^3u$. Thus the configuration of
Theorem~\ref{thm:T1} and its difference multiset are classical; the results of
Sections~\ref{sec:T2}--\ref{sec:T4} concern its maximality, its rigidity and the associated cosine sum.
\end{remark}

\begin{corollary}[golden ladder]\label{cor:ladder}
Let $K>0$, $u=\nu_{203}$ and $F>\nu_{200}$. Then
\[
  \Dl(S)=\{\nu_{200},\nu_{201},\nu_{202},\nu_{203}\},\qquad
  r_S(\nu_{200})=r_S(\nu_{203})=1,\quad r_S(\nu_{201})=r_S(\nu_{202})=2.
\] More generally, for $m\in\Z$ and
$u=\nu_{m+3}$ and $F>\nu_m$ one has $\Dl(S)=\{\nu_m,\nu_{m+1},\nu_{m+2},\nu_{m+3}\}$.
\end{corollary}

\begin{proof}
By Lemma~\ref{lem:ladder}, $\varphi^k\nu_{m+3}=\nu_{m+3-k}$, so
$\Dset(\nu_{m+3})=\{\nu_{m+3},\nu_{m+2},\nu_{m+1},\nu_m\}$; apply Theorem~\ref{thm:T1}. For $m=200$
note $a+b=\nu_{200}$, so $F>\nu_{200}$ is the hypothesis $F>a+b$.
\end{proof}

\begin{remark}[structure of $S$]\label{rem:structure}
\begin{enumerate}[label=(\alph*)]
\item As a set, $S=(F-a-b)+\bigl(\{0,a\}+\{0,b\}\bigr)=(F-\varphi^3u)+u\{0,\varphi,\varphi^2,\varphi^3\}$,
  by Lemma~\ref{lem:cube} and $\varphi+\varphi^2=\varphi^3$. Thus $S$ is a translate of the dilated
  two-dimensional cube $u\bigl(\{0,\varphi\}+\{0,\varphi^2\}\bigr)$. By Lemma~\ref{lem:cube}, for every
  two-dimensional cube with $\beta\neq2\alpha$ the multiplicities form the multiset $\{1,1,2,2\}$; listed
  in increasing order of the four values they are $1,2,2,1$ if $\alpha<\beta<2\alpha$ (order
  $\beta-\alpha<\alpha<\beta<\alpha+\beta$) and $2,1,2,1$ if $\beta>2\alpha$ (order
  $\alpha<\beta-\alpha<\beta<\alpha+\beta$). For the golden quadruple $\beta/\alpha=\varphi\in(1,2)$, which gives
  $1,2,2,1$; the golden ratio is what turns the four distinct values into a geometric progression.
\item Every three-element subset of $S$ has a three-term geometric progression with ratio $\varphi$ as
  its difference set: the differences of $\{s_1,s_2,s_3\}$, $\{s_2,s_3,s_4\}$, $\{s_1,s_2,s_4\}$ and
  $\{s_1,s_3,s_4\}$ are $\{u,\varphi u,\varphi^2u\}$, $\{u,\varphi u,\varphi^2u\}$,
  $\{\varphi u,\varphi^2u,\varphi^3u\}$ and $\{\varphi u,\varphi^2u,\varphi^3u\}$, respectively.
\item $S$ is not a Sidon set, since $r_S(\varphi u)=2$. Its full difference set has
  $|S-S|=2|\Dl(S)|+1=9$ elements; for every four-element $A\subset\R$ one has $7\le|A-A|\le13$, the lower
  bound being attained exactly by arithmetic progressions and the upper bound exactly by Sidon sets
  \cite{Nathanson1996,TaoVu2006}. So $S$ lies near the arithmetic end of this range.
\item The difference multiset does not depend on $F$; the condition $F>a+b$ only guarantees positivity.
\end{enumerate}
\end{remark}

\section{Maximality and uniqueness}\label{sec:T2}

\begin{theorem}\label{thm:T2}
Let $u>0$ and $\Dset=\{u,\varphi u,\varphi^2u,\varphi^3u\}$.
\begin{enumerate}[label=\textup{(\roman*)}]
\item There do not exist five (or more) distinct real numbers all of whose positive differences belong to
  $\Dset$.
\item If $A\subset\R$ has four elements and $\Dl(A)\subseteq\Dset$, then
  $A=\{y,\;y+\varphi u,\;y+\varphi^2u,\;y+\varphi^3u\}$ with $y=\min A$; that is, $A$ is a translate of the golden
  quadruple $S(F,u)$ (for any admissible $F$). In particular $\Dl(A)=\Dset$ and $r_A=r_S$.
\end{enumerate}
Hence $\max\{|A|:\ A\subseteq\R,\ \Dl(A)\subseteq\Dset\}=4$, attained exactly by the translates of $S$.
\end{theorem}

\begin{proof}
Let $y_1<\dots<y_m$ be distinct reals with all positive differences in $\Dset$, and let
$g_i=y_{i+1}-y_i$ ($1\le i<m$) be the gaps. Each $g_i\in\Dset$, so $g_i\ge u$, and the span
$y_m-y_1=g_1+\dots+g_{m-1}$ also lies in $\Dset$.

(i) Suppose $m\ge5$. The five numbers $y_1<\dots<y_5$ still have the property, so we may assume $m=5$.
Then $y_5-y_1\ge4u$. By Lemma~\ref{lem:phi}(iii), $\varphi^2u<3u<4u<\varphi^3u$, so the only element of
$\Dset$ that is $\ge4u$ is $\varphi^3u$; hence $g_1+g_2+g_3+g_4=\varphi^3u=(2\varphi+1)u$.
If all $g_i=u$, then $4=2\varphi+1$, i.e.\ $\varphi=3/2$, contradicting Lemma~\ref{lem:phi}(iii).
Otherwise some $g_i\ge\varphi u$, and the sum is at least $(3+\varphi)u>(2\varphi+1)u$, because $\varphi<2$.
In both cases we reach a contradiction.

(ii) Now $m=4$. The span is $\ge3u>\varphi^2u$, so $g_1+g_2+g_3=\varphi^3u$. No gap equals $\varphi^3u$
(the other two gaps are positive). If some gap equals $\varphi^2u$, the other two sum to
$(\varphi^3-\varphi^2)u=\varphi u<2u$, which is impossible since each is $\ge u$. Hence every gap lies in
$\{u,\varphi u\}$. If exactly $k$ of them equal $\varphi u$, then $(3-k)+k\varphi=1+2\varphi$, that is
$(2-k)(1-\varphi)=0$, so $k=2$. Thus $(g_1,g_2,g_3)$ is a permutation of $(u,\varphi u,\varphi u)$.
If $(g_1,g_2,g_3)=(u,\varphi u,\varphi u)$ then $y_4-y_2=2\varphi u$, and if
$(g_1,g_2,g_3)=(\varphi u,\varphi u,u)$ then $y_3-y_1=2\varphi u$; but $2\varphi u\notin\Dset$, since
$2\varphi\ne1$, $2\varphi\ne\varphi$, $2\varphi\ne\varphi+1$ (as $\varphi\ne1$) and $2\varphi\ne2\varphi+1$.
Therefore $(g_1,g_2,g_3)=(\varphi u,u,\varphi u)$ and
$A=\{y_1,\,y_1+\varphi u,\,y_1+(\varphi+1)u,\,y_1+(2\varphi+1)u\}=\{y_1,y_1+\varphi u,y_1+\varphi^2u,y_1+\varphi^3u\}$.
By Remark~\ref{rem:structure}(a), $A=S(F,u)+(y_1-F+\varphi^3u)$, and Theorem~\ref{thm:T1} gives
$\Dl(A)=\Dset$ and $r_A=r_S$. The final sentence combines (i), (ii) and Theorem~\ref{thm:T1}.
\end{proof}

\begin{remark}
A counting argument alone gives only $|A|\le5$: the differences $y_i-y_1$ ($2\le i\le m$) are $m-1$
distinct elements of the four-element set $\Dset$. The content of Theorem~\ref{thm:T2}(i) is that the
value $5$ is impossible.
\end{remark}

The same conclusion holds for the whole multiplicative orbit $\{c\varphi^n:n\in\Z\}$. The key is the
following ``golden unit equation''.

\begin{lemma}\label{lem:unit}
Let $p,q,r\in\Z$. Then $\varphi^p+\varphi^q=\varphi^r$ if and only if $\{p,q\}=\{r-1,r-2\}$. In particular
$p\neq q$ and $|p-q|=1$.
\end{lemma}

\begin{proof}
``If'' is Lemma~\ref{lem:phi}(ii). Conversely, assume $p\ge q$ and put $\alpha=p-r\ge\beta=q-r$, so that
$\varphi^\alpha+\varphi^\beta=1$. If $\alpha\ge0$ the left side exceeds $1$. If $\alpha\le-2$, then
$\varphi^\alpha+\varphi^\beta\le2\varphi^{-2}=2(2-\varphi)=4-2\varphi<1$, because $\varphi>3/2$. Hence
$\alpha=-1$ and $\varphi^\beta=1-\varphi^{-1}=2-\varphi=\varphi^{-2}$ (Lemma~\ref{lem:phi}(ii)), so $\beta=-2$
by Lemma~\ref{lem:phi}(vi).
\end{proof}

\begin{remark}\label{rem:nagell}
Lemma~\ref{lem:unit} is classical. Nagell \cite[Th\'eor\`eme~1]{Nagell1964} states that in a quadratic field a relation
$1+\xi+\xi_1=0$ in units exists only for $\Q(\sqrt{-3})$ and $\Q(\sqrt5)$, and in $\Q(\sqrt5)$ only as
$1+\xi-\xi^2=0$, $1-\xi+\xi^{-1}=0$, $1-\xi^{-1}-\xi^{-2}=0$, $\xi=\varphi$ or $\bar\varphi$ (with the proof
in an earlier paper cited there). Lemma~\ref{lem:unit} is equivalent to $\varphi^{p-r}+\varphi^{q-r}=1$, the case of positive units in
the unit equation $\xi+\xi'=1$ in $\Z[\varphi]^\times$. A unit $\xi$ such that $1-\xi$ is also a unit is
called \emph{exceptional}; Nagell determined the exceptional units of all number fields with unit rank at most
one (see \cite[\S1, no.~2]{Nagell1969b} and the works cited there, and \cite[\S3]{Lenstra1977}). In $\Q(\sqrt5)$ they are the six units
$\varphi,\ \varphi^2,\ \varphi^{-1},\ \varphi^{-2},\ -\varphi,\ -\varphi^{-1}$ (indeed $1-\varphi=-\varphi^{-1}$,
$1-\varphi^2=-\varphi$, $1-\varphi^{-1}=\varphi^{-2}$, $1-\varphi^{-2}=\varphi^{-1}$, $1+\varphi=\varphi^2$, $1+\varphi^{-1}=\varphi$). Conversely, writing $\xi=\pm\varphi^p$ and $1-\xi=\pm\varphi^q$, the relation
$\xi+(1-\xi)=1$ becomes one of $\varphi^p+\varphi^q=\varphi^0$, $\varphi^0+\varphi^q=\varphi^p$, $\varphi^0+\varphi^p=\varphi^q$ (the
fourth sign pattern, $-\varphi^p-\varphi^q=1$, is impossible), and Lemma~\ref{lem:unit} leaves exactly the six units listed.
The same list yields Theorem~2$'$(ii): after translating and dilating so that two of the
points are $0$ and $1$ (Lenstra's normalisation \cite[\S2, p.~243]{Lenstra1977}; cf.\ \cite{LeutbecherMartinet1982}), each remaining point
$\xi$ is an exceptional unit, and a finite check over the six units above gives the claim.
We have not found this stated explicitly; our proof is direct.
\end{remark}

\begin{theoremTwoPrime}
Let $c>0$ and $\Phi_c=\{c\varphi^n:n\in\Z\}$.
\begin{enumerate}[label=\textup{(\roman*)}]
\item A three-element set $\{y_1<y_2<y_3\}$ has $\Dl(\{y_1,y_2,y_3\})\subseteq\Phi_c$ if and only if its gaps are
  $(c\varphi^{n+1},c\varphi^n)$ or $(c\varphi^n,c\varphi^{n+1})$ for some $n\in\Z$.
\item A four-element set $A$ has $\Dl(A)\subseteq\Phi_c$ if and only if
  $A=t+c\varphi^n\{0,\varphi,\varphi^2,\varphi^3\}$ for some $t\in\R$ and $n\in\Z$, i.e.\ $A$ is a translate of a
  golden quadruple $S(F,c\varphi^n)$.
\item No set of five or more distinct reals has all its positive differences in $\Phi_c$.
\end{enumerate}
In particular Theorem~\textup{\ref{thm:T2}} follows again: if $\Dl(A)\subseteq\Dset(u)\subseteq\Phi_u$ and
$|A|=4$, then $A=t+u\varphi^n\{0,\varphi,\varphi^2,\varphi^3\}$ and $\Dl(A)=\{u\varphi^n,\dots,u\varphi^{n+3}\}\subseteq\Dset(u)$
forces $n=0$.
\end{theoremTwoPrime}

\begin{proof}
Write the gaps of an increasing set as $g_i=c\varphi^{e_i}$.

(i) The differences are $g_1,g_2,g_1+g_2$, so we need $\varphi^{e_1}+\varphi^{e_2}=\varphi^r$ for some $r$.
By Lemma~\ref{lem:unit} this means $\{e_1,e_2\}=\{r-1,r-2\}$, which is the claim with $n=r-2$; conversely
such gaps give $g_1+g_2=c\varphi^{n+2}$.

(ii) Let $A=\{y_1<y_2<y_3<y_4\}$ with $\Dl(A)\subseteq\Phi_c$, and write $g_1+g_2=c\varphi^r$. Applying (i) to the
triples $\{y_1,y_2,y_3\}$, $\{y_1,y_3,y_4\}$ and $\{y_2,y_3,y_4\}$ gives, respectively,
$\{e_1,e_2\}=\{r-1,r-2\}$, $e_3\in\{r-1,r+1\}$ and $|e_2-e_3|=1$. If $e_3=r+1$ then $|e_2-e_3|\ge2$, a
contradiction; so $e_3=r-1$, then $|e_2-(r-1)|=1$ with $e_2\in\{r-1,r-2\}$ forces $e_2=r-2$, and $e_1=r-1$.
With $n=r-2$ the gaps are $(c\varphi^{n+1},c\varphi^n,c\varphi^{n+1})$, so
$A=y_1+c\varphi^n\{0,\varphi,\varphi+1,2\varphi+1\}=y_1+c\varphi^n\{0,\varphi,\varphi^2,\varphi^3\}$. Conversely,
such a set has $\Dl(A)=\{c\varphi^n,\dots,c\varphi^{n+3}\}\subseteq\Phi_c$ by Theorem~\ref{thm:T1} with
$u=c\varphi^n$, and it is a translate of $S(F,c\varphi^n)$ by Remark~\ref{rem:structure}(a).

(iii) Let $y_1<\dots<y_5$ have all differences in $\Phi_c$. By (ii) applied to $\{y_1,\dots,y_4\}$ there is
$s$ with $(e_1,e_2,e_3)=(s+1,s,s+1)$, and applied to $\{y_2,\dots,y_5\}$ there is $s'$ with
$(e_2,e_3,e_4)=(s'+1,s',s'+1)$. Then $e_2=s=s'+1$ and $e_3=s+1=s'$, whence $s'=s+1=s'+2$, a contradiction.
Sets with more than five elements contain five-element subsets with the same property.
\end{proof}

Part~(i) of Theorem~2$'$ is a restatement of Lemma~\ref{lem:unit} and hence classical (Remark~\ref{rem:nagell}).
Part~(iii) is known, see
Remark~\ref{rem:lenstra}.

\begin{definition}\label{def:N}
For $x>1$ let $N(x)\in\{2,3,\dots\}\cup\{\infty\}$ be the supremum of the cardinalities of finite sets $A\subset\R$ all of
whose positive differences lie in $\{cx^n:n\in\Z\}$ for some $c>0$.
\end{definition}

\begin{remark}[relation to Lenstra's constant]\label{rem:lenstra}
The bound in Theorem~2$'$(iii), and hence in Theorem~\ref{thm:T2}(i), i.e.\ $N(\varphi)\le4$, is known, and so is
$N(\varrho)\le5$ for the plastic number $\varrho$, $\varrho^3=\varrho+1$ (Section~\ref{sec:literature}). We recall the
argument, which is that of Lenstra \cite[(2.2)]{Lenstra1977}. Let $A$ be finite with all positive differences in
$\{cx^n:n\in\Z\}$, $x\in\{\varphi,\varrho\}$, and put $A'=c^{-1}(A-\min A)$. Every element of $A'$ is a sum of gaps of
the form $x^e$, $e\in\Z$. Since $\varphi^{-1}=\varphi-1$ and $\varrho^{-1}=\varrho^2-1$, all integral powers of $x$ lie
in $\Z[x]$; hence $A'\subset\Z[x]$, and all differences of elements of $A'$ are units $\pm x^n$ of $\Z[x]$. Any ring
homomorphism $\psi$ from $\Z[x]$ to a finite field maps units to non-zero elements, so $\psi$ is injective on $A'$
and $|A'|$ is at most the size of the field. Since $X^2-X-1\equiv X^2+X+1\pmod2$, the assignment $\varphi\mapsto\omega$
defines a surjective ring homomorphism $\Z[\varphi]\cong\Z[X]/(X^2-X-1)\to\mathbb F_4=\mathbb F_2[\omega]$,
$\omega^2+\omega+1=0$; hence $N(\varphi)\le4$. Since $2^3-2-1=5$, the evaluation $\varrho\mapsto2$ defines a surjective
ring homomorphism $\Z[\varrho]\cong\Z[X]/(X^3-X-1)\to\mathbb F_5$ (indeed $X^3-X-1\equiv(X-2)(X^2+2X+3)\pmod5$, the
quadratic factor being irreducible because its discriminant $-8\equiv2$ is not a square modulo~$5$); hence
$N(\varrho)\le5$. (In the notation of \cite[(2.1)]{Lenstra1977}, $L=4$ for $\Q(\sqrt5)$ and $L=5$ for $\Q(\varrho)$; these are the bounds $M\le L$ of \cite[(2.2)]{Lenstra1977}.) Together with Theorem~\ref{thm:T1} and Proposition~\ref{prop:plastic} this gives $N(\varphi)=4$
and $N(\varrho)=5$. In Lenstra's notation these are the values $M=4$ for $\Q(\sqrt5)$ \cite[(3.9)]{Lenstra1977}; cf.\ \cite[\S3.0]{LeutbecherMartinet1982},
obtained there by applying \cite[(2.4)(b)]{Lenstra1977} with $f=X^2-X-1$, i.e.\ from the sequence $0,1,\varphi,\varphi+1=\varphi^2$, which is
$\varphi^{-1}\{0,\varphi,\varphi^2,\varphi^3\}$, and $M=5$ for the complex cubic field $K=\Q(\alpha)$, $\alpha^3-\alpha-1=0$,
$\Delta=-23$ \cite[(3.10)]{Lenstra1977}, see also \cite[\S3.0, Thm.~4.1.1]{LeutbecherMartinet1982}, i.e.\ $K=\Q(\varrho)$; since $\operatorname{disc}(X^3-X-1)=-23$ is squarefree,
$\Z[\varrho]$ is its ring of integers. What is not contained
in~\cite{Lenstra1977} is the uniqueness statement of Theorem~\ref{thm:T2}(ii), the explicit uniqueness statement
of Theorem~2$'$(ii) (see above), the rigidity Theorem~3$'$, and the analysis of other bases $x$
(Section~\ref{sec:T3}).
\end{remark}

\section{Uniqueness of the base and closure residues}\label{sec:T3}

\begin{definition}\label{def:Sx}
Let $x>1$, $u>0$ and $F>(x+x^2)u$. The \emph{$x$-quadruple} is
\[
  S_x=\bigl(F,\;F-xu,\;F-x^2u,\;F-(x+x^2)u\bigr),
\]
and its \emph{target} is $\Dset_x=\{u,xu,x^2u,x^3u\}$. Put $P(x)=x^2-x-1$ and define the two
\emph{closure residues}
\[
  \rho_-(x)=(x^2u-xu)-u=P(x)\,u,\qquad \rho_+(x)=x^3u-(xu+x^2u)=x\,P(x)\,u .
\]
The residue $\rho_-$ measures the failure of the lower closure $b-a=u$, and $\rho_+$ the failure of
the upper closure $a+b=x^3u$, for $a=xu$, $b=x^2u$. By Lemma~\ref{lem:phi}(v), $S_\varphi=S$.
\end{definition}

\begin{theorem}\label{thm:T3}
Let $x>1$.
\begin{enumerate}[label=\textup{(\alph*)}]
\item The difference multiset of $S_x$ is
  \[
    \{(x^2-x)u,\ xu,\ xu,\ x^2u,\ x^2u,\ (x^2+x)u\}
    =\{u+\rho_-(x),\ xu,\ xu,\ x^2u,\ x^2u,\ x^3u-\rho_+(x)\}.
  \]
\item $\Dl(S_x)=\Dset_x$ $\iff$ $\rho_-(x)=0$ $\iff$ $\rho_+(x)=0$ $\iff$ $x=\varphi$.
\item If $x\ne\varphi$, then $\rho_-(x)\ne0$ and $\rho_+(x)\ne0$; both have the sign of $x-\varphi$,
  $\rho_+(x)=x\rho_-(x)$, and $|\rho_-(x)|=|x-\varphi|\,(x-\bphi)\,u$.
\item More strongly, $\Dl(S_x)$ is a geometric progression of four terms (with arbitrary first term and
  arbitrary ratio) if and only if $x=\varphi$. In particular, for $x\ne\varphi$ there are no $c>0$ and
  $k\in\Z$ with $\Dl(S_x)=\{cx^k,cx^{k+1},cx^{k+2},cx^{k+3}\}$.
\end{enumerate}
\end{theorem}

\begin{proof}
(a) Apply Lemma~\ref{lem:cube} with $\alpha=xu<\beta=x^2u$ to $S_x=F-Q$; then $\beta-\alpha=(x^2-x)u$ and
$\alpha+\beta=(x^2+x)u$. The second form is the definition of $\rho_\pm$.

(b) The last three conditions are equivalent because $\rho_-=P(x)u$, $\rho_+=xP(x)u$ with $x,u>0$, and
$P(x)=0\iff x=\varphi$ for $x>1$ (Lemma~\ref{lem:phi}(v)). If $x=\varphi$, Theorem~\ref{thm:T1} gives
$\Dl(S_\varphi)=\Dset_\varphi$. Conversely assume $\Dl(S_x)=\Dset_x$. If $x=2$ then $(x^2-x)u=xu$, so
$|\Dl(S_x)|=3<4=|\Dset_x|$, impossible. Let $x\ne2$. Then $(x^2-x)u$, $xu$, $x^2u$, $(x^2+x)u$ are four
distinct numbers ($x^2-x\ne x$ as $x\ne2$; $x^2-x<x^2<x^2+x$; $x<x^2$), and none of $(x^2\mp x)u$ lies in
$\{xu,x^2u\}$, while none of $u,x^3u$ lies in $\{xu,x^2u\}$. Removing the common elements $xu,x^2u$ from
both sides yields $\{(x^2-x)u,(x^2+x)u\}=\{u,x^3u\}$. Comparing the smaller and the larger elements,
$x^2-x=1$ and $x^2+x=x^3$; either equation is $P(x)=0$, so $x=\varphi$.

(c) This is Lemma~\ref{lem:phi}(v) together with $\rho_-=P(x)u$, $\rho_+=x\rho_-$.

(d) ``If'' is Theorem~\ref{thm:T1}. Suppose $\Dl(S_x)$ is a four-term geometric progression; then
$|\Dl(S_x)|=4$, so $x\ne2$.

\emph{Case $1<x<2$.} Here $x^2-x<x$, so the increasing order is
$(x^2-x)u<xu<x^2u<(x^2+x)u$, with consecutive ratios $\frac{1}{x-1}$, $x$, $\frac{x+1}{x}$. Equality of the
first two gives $x(x-1)=1$, i.e.\ $P(x)=0$, so $x=\varphi$ (and then indeed $\frac{x+1}{x}=x$).

\emph{Case $x>2$.} Here $x<x^2-x$, so the increasing order is $xu<(x^2-x)u<x^2u<(x^2+x)u$, with ratios
$x-1$, $\frac{x}{x-1}$, $\frac{x+1}{x}$. Equality of the first two gives $(x-1)^2=x$, i.e.\ $x^2-3x+1=0$,
whose only root exceeding $2$ is $(3+\sqrt5)/2=\varphi^2$. Equality of the first and the third gives
$x(x-1)=x+1$, i.e.\ $x^2-2x-1=0$, whose only root exceeding $2$ is $1+\sqrt2$. If $\varphi^2=1+\sqrt2$, then
$1+\sqrt5=2\sqrt2$, and squaring, $6+2\sqrt5=8$, i.e.\ $\sqrt5=1$, which is false. So no $x>2$ works.

The final assertion holds because $\{cx^k,\dots,cx^{k+3}\}$ is a four-term progression with ratio $x$.
\end{proof}

\subsection{Fibonacci quotients}

\begin{lemma}[Cassini's identity]\label{lem:cassini}
$\Fib_{n+1}\Fib_{n-1}-\Fib_n^2=(-1)^n$ for all $n\ge1$.
\end{lemma}

\begin{proof}
For $n=1$: $\Fib_2\Fib_0-\Fib_1^2=-1$. If $c_n=\Fib_{n+1}\Fib_{n-1}-\Fib_n^2$, then
$c_{n+1}=\Fib_{n+2}\Fib_n-\Fib_{n+1}^2=(\Fib_{n+1}+\Fib_n)\Fib_n-\Fib_{n+1}(\Fib_n+\Fib_{n-1})
=\Fib_n^2-\Fib_{n+1}\Fib_{n-1}=-c_n$.
\end{proof}

\begin{proposition}\label{prop:fibres}
For $n\ge2$ let $x_n=\Fib_{n+1}/\Fib_n$. Then $x_n>1$, $P(x_n)=(-1)^n/\Fib_n^2$, and therefore
\[
  \rho_-(x_n)=\frac{(-1)^n}{\Fib_n^2}\,u,\qquad \rho_+(x_n)=\frac{(-1)^n\Fib_{n+1}}{\Fib_n^3}\,u .
\]
In particular $x_n\neq\varphi$ for all $n$, the residues alternate in sign, and
$|\rho_-(x_n)|=u/\Fib_n^2\to0$.
\end{proposition}

\begin{proof}
For $n\ge2$, $\Fib_{n-1}\ge1$, so $\Fib_{n+1}=\Fib_n+\Fib_{n-1}>\Fib_n$. Next,
$\Fib_n^2P(x_n)=\Fib_{n+1}^2-\Fib_{n+1}\Fib_n-\Fib_n^2=\Fib_{n+1}(\Fib_{n+1}-\Fib_n)-\Fib_n^2
=\Fib_{n+1}\Fib_{n-1}-\Fib_n^2=(-1)^n$ by Lemma~\ref{lem:cassini}. The rest follows from
Definition~\ref{def:Sx}.
\end{proof}

\begin{remark}[the integer shadow]\label{rem:shadow}
For $n\ge3$ and an integer $N>\Fib_{n+2}$, Lemma~\ref{lem:cube} with $\alpha=\Fib_n$, $\beta=\Fib_{n+1}$
shows that the integer quadruple $(N,N-\Fib_n,N-\Fib_{n+1},N-\Fib_{n+2})$ has differences
$\Fib_{n-1},\Fib_n,\Fib_n,\Fib_{n+1},\Fib_{n+1},\Fib_{n+2}$: four consecutive Fibonacci numbers with the same
multiplicities $1,2,2,1$ (they are distinct because $n\ge3$). For instance $(N,N-2,N-3,N-5)$ has
differences $1,2,2,3,3,5$. Here the additive closures $\beta-\alpha=\Fib_{n-1}$ and $\alpha+\beta=\Fib_{n+2}$
are exact, but the four values do not form a geometric progression. A geometric progression $(cx^k)$
with $c>0$, $x>1$ satisfies the Fibonacci recursion $cx^{k+2}=cx^{k+1}+cx^k$ only for $x=\varphi$; the golden
quadruple is the unique case in which the Fibonacci-type additive relations and geometric scaling hold
at the same time.
\end{remark}

\subsection{Rigidity for arbitrary four-point sets}

Theorem~\ref{thm:T3} concerns quadruples of the special form $S_x$. The next result removes this
restriction.

\begin{theoremThreePrime}
Let $x>1$ and $c>0$. There exists a four-element set $T\subset\R$ with
\[
  \Dl(T)=V:=\{c,\;cx,\;cx^2,\;cx^3\}
\]
if and only if $x=\varphi$. In that case $T=t+c\{0,\varphi,\varphi^2,\varphi^3\}$ for some $t\in\R$; thus $T$ is a
translate of the golden quadruple $S(F,c)$, and $r_T(c)=r_T(c\varphi^3)=1$, $r_T(c\varphi)=r_T(c\varphi^2)=2$.
\end{theoremThreePrime}

\begin{proof}
``If'' is Theorem~\ref{thm:T1} with $u=c$ and Remark~\ref{rem:structure}(a). Conversely, let
$T=\{y_1<y_2<y_3<y_4\}$ with $\Dl(T)=V$, gaps $g_1,g_2,g_3$, and put $\sigma_1=g_1+g_2=y_3-y_1$,
$\sigma_2=g_2+g_3=y_4-y_2$. The span $y_4-y_1$ is the largest difference, so it equals $\max V=cx^3$. All of
$g_1,g_2,g_3,\sigma_1,\sigma_2$ lie in $V$; each gap is $\ge\min V=c$, so $\sigma_1,\sigma_2\ge2c>c$; and $\sigma_1,\sigma_2<cx^3$. Hence
$\sigma_1,\sigma_2\in\{cx,cx^2\}$.

\emph{Step 1: $\sigma_1=\sigma_2=cx^2$.} Suppose $\sigma_1=cx$. Then $g_1,g_2\in V$ and $g_1,g_2<\sigma_1=cx$; the only element of $V$
below $cx$ is $c$, so $g_1=g_2=c$ and $x=2$. But then $g_3=cx^3-\sigma_1=6c\notin\{c,2c,4c,8c\}=V$, a
contradiction. The case $\sigma_2=cx$ is symmetric ($g_2=g_3=c$, $x=2$, $g_1=6c\notin V$).

\emph{Step 2.} Now $g_3=cx^3-\sigma_1=cx^2(x-1)$, $g_1=cx^3-\sigma_2=cx^2(x-1)$ and $g_2=\sigma_1-g_1=cx^2(2-x)$. Since $g_2>0$,
$x<2$. Since $g_3\in V$, $x^2(x-1)\in\{1,x,x^2,x^3\}$. The value $x^3$ is impossible, and $x^2$ would give
$x=2$. The value $x$ gives $x^2-x-1=0$, i.e.\ $x=\varphi$. The value $1$ gives $x^3-x^2=1$.

\emph{Step 3: $x^3-x^2=1$ is impossible.} In that case $g_1=g_3=c$ and $g_2=c(2x^2-x^3)=c(x^2-1)$, and
$g_2\in V$ requires $x^2-1\in\{1,x,x^2,x^3\}$. Now $x^2-1=1$ gives $x=\sqrt2$, but
$(\sqrt2)^3-(\sqrt2)^2=2\sqrt2-2\ne1$ (as $8\ne9$); $x^2-1=x$ gives $x=\varphi$, but
$\varphi^3-\varphi^2=\varphi\ne1$; $x^2-1=x^2$ is absurd; and $x^2-1=x^3$ gives $x^3-x^2=-1\ne1$.

\emph{Step 4: $x=\varphi$.} By Lemma~\ref{lem:phi}(ii), $g_1=g_3=c\varphi^2(\varphi-1)=c\varphi$ and
$g_2=c\varphi^2(2-\varphi)=c$, so $T=\{y_1,\,y_1+c\varphi,\,y_1+c\varphi^2,\,y_1+c\varphi^3\}$. The
multiplicities follow from Theorem~\ref{thm:T1}.
\end{proof}

\begin{remark}
The proof shows that apart from $x=\varphi$, the only ``candidates'' surviving the closure of the span are $x=2$ and the real
root $x\approx1.465\,57$ of $x^3=x^2+1$. The value $x=2$ fails because either $g_3=6c\notin V$
(Step~1) or $g_2=cx^2(2-x)=0$ (Step~2); the root of $x^3=x^2+1$ fails because the middle gap
$g_2=c(x^2-1)$ does not lie in $V$ (Step~3). For $0<x<1$ a progression with
ratio $x$ is a progression with ratio $1/x$ read backwards, so for ratios $x\ne1$ the answer is
$x\in\{\varphi,\varphi^{-1}\}$. For three points the analogue is immediate: if $\{y_1<y_2<y_3\}$ has
$\Dl=\{c,cx,cx^2\}$ with $x>1$, then the span is $cx^2$, the two gaps are distinct (otherwise there are
only two differences), so $\{g_1,g_2\}=\{c,cx\}$ and $c+cx=cx^2$, i.e.\ $x=\varphi$.
\end{remark}

\begin{remark}[Relation to {\cite[Lemma~3]{CampolinaMailybaev2021}}]\label{rem:CM}
The proof of Theorem~3$'$ is self-contained and does not use \cite{CampolinaMailybaev2021}. Each three-point
subset of a set $T$ as in Theorem~3$'$ gives a relation $x^l=x^k+1$, $l>k\ge0$, i.e.\ $x$
is a root of $X^l-X^k-1$. By \cite[Lemma~3]{CampolinaMailybaev2021} (with a computer-assisted step), the
common roots $x\ge1.05$ of two distinct such trinomials are exactly $\varrho^{1/m}$,
$m=1,\dots,5$. For $\varphi$ no such pair exists: by Lemma~\ref{lem:unit}, $X^2-X-1$ is the only
trinomial $X^l-X^k-1$ ($l>k\ge0$) vanishing at $\varphi$, and all four three-point
relations in $\{0,\varphi,\varphi^2,\varphi^3\}$ reduce to it. This is consistent with
\cite[Lemma~3]{CampolinaMailybaev2021} and needs no computation. Theorem~3$'$ can also be derived from that lemma by a short argument. The span $cx^3$ of $T$ is a sum of three gaps, each
$\ge c$, so $x\ge\sqrt[3]3>\varrho$; for such $x$ at most one trinomial $X^l-X^k-1$ vanishes at $x$ (this step is even elementary: as $x\ge\sqrt[3]3>\sqrt2$ (because $3^2>2^3$), $x^{l-k}=1+x^{-k}\le2<x^2$ forces $l=k+1$, and $x^k(x-1)=1$ determines $k$ because $x^k(x-1)$ is strictly increasing in $k$; thus \cite[Lemma~3]{CampolinaMailybaev2021} is used only for the uniqueness of the trinomial, which here is elementary), so all three-point
subsets of $T$ are similar up to reflection. Normalising $T=\{0,p,q,1\}$, the triples $\{0,p,1\}$ and $\{0,q,1\}$ give
$q=1-p$ with $0<p<\frac12$, and the triple $\{0,p,1-p\}$, with parts $p/(1-p)$ and $(1-2p)/(1-p)$, forces
$p/(1-p)=1-p$ (as $p/(1-p)=p$ would give $p=0$), i.e.\ $p=(1-p)^2$, $p=\varphi^{-2}$; the distances are then $\varphi^{-3},\varphi^{-2},\varphi^{-1},1$, so
$x=\varphi$.
\end{remark}

\section{The cosine sum: non-periodic, but almost periodic}\label{sec:T4}

\begin{lemma}\label{lem:square}
For finite $A\subset\R$ and $t\in\R$: $|\widehat A(t)|^2=|A|+2C_A(t)$. Consequently
$-|A|/2\le C_A(t)\le\binom{|A|}{2}$, and $C_A(t)=\binom{|A|}{2}$ if and only if $dt\in\Z$ for every
$d\in\Dl(A)$.
\end{lemma}

\begin{proof}
$|\widehat A(t)|^2=\sum_{x,y\in A}e((x-y)t)=|A|+\sum_{x>y}\bigl(e((x-y)t)+e(-(x-y)t)\bigr)=|A|+2C_A(t)$.
The lower bound follows from $|\widehat A|^2\ge0$. The upper bound holds because $C_A$ is a sum of
$\binom{|A|}2$ cosines, and equality holds if and only if each cosine equals $1$, i.e.\ $(x-y)t\in\Z$ for
all $x>y$ in $A$.
\end{proof}

\begin{proposition}[closed form]\label{prop:closed}
For all $t\in\R$,
\[
  C_S(t)=\cos(2\pi ut)+2\cos(2\pi\varphi ut)+2\cos(2\pi\varphi^2ut)+\cos(2\pi\varphi^3ut)
        =8\cos^2(\pi at)\cos^2(\pi bt)-2 .
\]
Equivalently, for the exponential sum $\widehat S(t)=\sum_{x\in S}e(xt)$ of the four points of $S$,
$|\widehat S(t)|^2=4+2C_S(t)=16\cos^2(\pi at)\cos^2(\pi bt)$, i.e.\ $|\widehat S(t)|=4\,|\cos(\pi at)\cos(\pi bt)|$.
Consequently:
\begin{enumerate}[label=\textup{(\roman*)}]
\item $-2\le C_S(t)\le6$ for all $t$;
\item $C_S(t)=6$ if and only if $t=0$;
\item $C_S(t)=-2$ if and only if $t\in\frac{2\Z+1}{2a}\cup\frac{2\Z+1}{2b}$.
\end{enumerate}
\end{proposition}

\begin{proof}
The first expression is Definition~\ref{def:cosine} combined with Theorem~\ref{thm:T1}. Since
$S=F-(\{0,a\}+\{0,b\})$,
\[
  \widehat S(t)=e(Ft)\sum_{q\in\{0,a,b,a+b\}}e(-qt)=e(Ft)\,\bigl(1+e(-at)\bigr)\bigl(1+e(-bt)\bigr).
\]
As $|1+e(\theta)|^2=2+2\cos2\pi\theta=4\cos^2\pi\theta$, we get $|\widehat S(t)|^2=16\cos^2(\pi at)\cos^2(\pi bt)$,
and Lemma~\ref{lem:square} with $|S|=4$ gives $C_S=(|\widehat S|^2-4)/2$, which is the closed form.
(i) follows from $0\le\cos^2\le1$.
(ii) $C_S(t)=6$ if and only if $\cos^2(\pi at)=\cos^2(\pi bt)=1$, i.e.\ $at\in\Z$ and $bt\in\Z$. If $t\ne0$,
then $at\ne0$ and $b/a=(bt)/(at)\in\Q$, contradicting $b/a=\varphi\notin\Q$ (Lemma~\ref{lem:phi}(iv)).
(iii) $C_S(t)=-2$ if and only if $\cos(\pi at)=0$ or $\cos(\pi bt)=0$, i.e.\ $at\in\frac12+\Z$ or
$bt\in\frac12+\Z$.
\end{proof}

\begin{lemma}[trigonometric polynomials are almost periodic]\label{lem:trigAP}
Let $\lambda_1,\dots,\lambda_k\in\R$, $c_1,\dots,c_k\in\C$ and $P(t)=\sum_{j=1}^kc_je(\lambda_jt)$.
\begin{enumerate}[label=\textup{(\alph*)}]
\item For all $t,\tau\in\R$: $|P(t+\tau)-P(t)|\le2\pi\sum_j|c_j|\,\nint{\lambda_j\tau}$. Likewise
  $|\cos(2\pi\lambda(t+\tau))-\cos(2\pi\lambda t)|\le2\pi\nint{\lambda\tau}$ for real $\lambda$.
\item For every $\delta>0$ the set $B(\delta)=\{\tau\in\R:\ \nint{\lambda_j\tau}<\delta\ \text{for all } j\}$ is
  relatively dense.
\item $P$ is almost periodic in the sense of Bohr; so is every finite real linear combination of functions
  $t\mapsto\cos(2\pi\lambda t)$.
\end{enumerate}
\end{lemma}

\begin{proof}
(a) $P(t+\tau)-P(t)=\sum_jc_je(\lambda_jt)\bigl(e(\lambda_j\tau)-1\bigr)$; apply the triangle inequality and
Lemma~\ref{lem:nint}. For the cosine, $|\cos(2\pi\lambda(t+\tau))-\cos(2\pi\lambda t)|
=|\operatorname{Re}\bigl(e(\lambda t)(e(\lambda\tau)-1)\bigr)|\le|e(\lambda\tau)-1|$.

(b) Let $\T^k=\R^k/\Z^k$ with quotient map $\pi$, and define
$d(\pi(\xi),\pi(\eta))=\max_j\nint{\xi_j-\eta_j}$. By Lemma~\ref{lem:nint} this is a well-defined,
translation-invariant metric on $\T^k$, and $d(\pi(\xi),\pi(\eta))\le\max_j|\xi_j-\eta_j|$, so $\pi$ is
continuous and $\T^k=\pi([0,1]^k)$ is compact. The map $\Theta:\R\to\T^k$,
$\Theta(s)=\pi(\lambda_1s,\dots,\lambda_ks)$, is an additive homomorphism. Let $\mathcal K$ be the closure of
$\Theta(\R)$; it is compact. The open balls of radius $\delta$ centred at the points $\Theta(s)$, $s\in\R$, cover
$\mathcal K$, so finitely many of them, centred at $\Theta(s_1),\dots,\Theta(s_r)$, already cover $\mathcal K$.
Put $R=\max_i|s_i|$ and let $\sigma\in\R$. Since $\Theta(-\sigma)\in\mathcal K$, there is $i$ with
$d(\Theta(-\sigma),\Theta(s_i))<\delta$; adding $\Theta(\sigma)$ and using translation invariance,
$d(0,\Theta(\sigma+s_i))<\delta$, i.e.\ $\tau=\sigma+s_i\in B(\delta)$, and $|\tau-\sigma|\le R$. Hence every closed
interval of length $L=2R+1$ meets $B(\delta)$.

(c) $P$ is continuous. Given $\varepsilon>0$, choose $\delta=\varepsilon/\bigl(2\pi(1+\sum_j|c_j|)\bigr)$; by (a),
$B(\delta)\subseteq E(P;\varepsilon)$, and a set containing a relatively dense set is relatively dense. The last
claim follows from $\cos(2\pi\lambda t)=\frac12e(\lambda t)+\frac12e(-\lambda t)$.
\end{proof}

\begin{proposition}[periodicity criterion]\label{prop:periodic}
Let $A\subset\R$ be finite with $|A|\ge2$. The cosine sum $C_A(t)=\sum_{d\in\Dl(A)}r_A(d)\cos(2\pi dt)$ is
periodic if and only if all ratios of its exponents $d$ are rational, i.e.\ $d/d'\in\Q$ for all
$d,d'\in\Dl(A)$.
\end{proposition}

This criterion is well known. It is a special case of a theorem of Mirotin and Mirotin
\cite[\S2]{MirotinMirotin2009}: a sum of continuous, nonconstant periodic functions on $\R$, none of whose proper
subsums is constant, is periodic if and only if the periods of the summands are commensurable. Here the summands
$r_A(d)\cos(2\pi dt)$ have periods $1/d$, and $1/d$, $1/d'$ are commensurable exactly when $d/d'\in\Q$; moreover
every nonempty subsum $\sum_{d\in E}r_A(d)\cos(2\pi dt)$, $\emptyset\ne E\subseteq\Dl(A)$, is nonconstant, because
its value at $t=0$ is $\sum_{d\in E}r_A(d)>0$ while its mean $\frac1T\int_0^T$ tends to $0$ as $T\to\infty$. For
completeness we include the short direct proof.

\begin{proof}
($\Leftarrow$) Fix $d_0\in\Dl(A)$ and write $d=(p_d/q_d)\,d_0$ with $p_d,q_d$ positive integers. With
$Q=\operatorname{lcm}_d q_d$ and $T=Q/d_0$ we get $dT=p_dQ/q_d\in\Z$, so every term $\cos(2\pi dt)$ has period $T$.
($\Rightarrow$) If $T>0$ is a period, then $C_A(T)=C_A(0)=\binom{|A|}2$, so $dT\in\Z\setminus\{0\}$ for every
$d\in\Dl(A)$ by Lemma~\ref{lem:square}, and $d/d'=(dT)/(d'T)\in\Q$.
\end{proof}

\begin{theorem}\label{thm:T4}
Let $C_S(t)=\sum_{1\le i<j\le4}\cos\bigl(2\pi(s_i-s_j)t\bigr)
=\cos(2\pi ut)+2\cos(2\pi\varphi ut)+2\cos(2\pi\varphi^2ut)+\cos(2\pi\varphi^3ut)$.
\begin{enumerate}[label=\textup{(\alph*)}]
\item $C_S$ is not periodic.
\item $C_S$ is almost periodic in the sense of Bohr. Quantitatively, for all $t,\tau\in\R$,
\[
  |C_S(t+\tau)-C_S(t)|\le2\pi\bigl(4\nint{u\tau}+6\nint{\varphi u\tau}\bigr),
\]
  so $E(C_S;\varepsilon)$ contains the relatively dense set
  $\{\tau:\ \nint{u\tau}<\varepsilon/(20\pi),\ \nint{\varphi u\tau}<\varepsilon/(20\pi)\}$.
\end{enumerate}
\end{theorem}

\begin{proof}
(a) If $T>0$ were a period, then $C_S(T)=C_S(0)=6$, contradicting Proposition~\ref{prop:closed}(ii).
(Alternatively: $u$ and $\varphi u$ are incommensurable, and Proposition~\ref{prop:periodic} applies.)

(b) By Lemma~\ref{lem:trigAP}(a), $|C_S(t+\tau)-C_S(t)|\le2\pi\sum_{d}r_S(d)\nint{d\tau}$. Put
$\theta_1=u\tau$, $\theta_2=\varphi u\tau$. Then $\varphi^2u\tau=\theta_1+\theta_2$ and $\varphi^3u\tau=\theta_1+2\theta_2$, so
by Lemma~\ref{lem:nint}
\[
  \sum_dr_S(d)\nint{d\tau}\le\nint{\theta_1}+2\nint{\theta_2}+2\bigl(\nint{\theta_1}+\nint{\theta_2}\bigr)
  +\bigl(\nint{\theta_1}+2\nint{\theta_2}\bigr)=4\nint{\theta_1}+6\nint{\theta_2}.
\]
If $\nint{\theta_1},\nint{\theta_2}<\varepsilon/(20\pi)$, the right side of the displayed estimate is $<\varepsilon$. The set
of such $\tau$ is $B(\varepsilon/(20\pi))$ for $(\lambda_1,\lambda_2)=(u,\varphi u)$, which is relatively dense by
Lemma~\ref{lem:trigAP}(b). Since $C_S$ is continuous, it is almost periodic.
\end{proof}

\begin{proposition}[explicit almost periods]\label{prop:fibAP}
For $n\in\N$ let $\tau_n=\Fib_n/u$. Then:
\begin{enumerate}[label=\textup{(\roman*)}]
\item $\Fib_{n+1}-\varphi\Fib_n=(-\varphi^{-1})^n$, hence $\nint{\varphi\Fib_n}\le\varphi^{-n}$, with equality
  for $n\ge2$ (for $n=1$, $\nint{\varphi}=\varphi^{-2}$);
\item $\sup_{t\in\R}|C_S(t+\tau_n)-C_S(t)|\le12\pi\varphi^{-n}$;
\item $C_S(\tau_n)=1+4\cos(2\pi\varphi^{-n})+\cos(4\pi\varphi^{-n})$, so that $C_S(\tau_n)<6$ for $n\ge1$ and
  $0<6-C_S(\tau_n)\le16\pi^2\varphi^{-2n}$.
\end{enumerate}
Thus the maximum value $6$, attained only at $t=0$, is approached along the explicit sequence
$\tau_n\to\infty$.
\end{proposition}

\begin{proof}
(i) Let $\epsilon_n=\Fib_{n+1}-\varphi\Fib_n$. Then $\epsilon_0=1$ and, using $1-\varphi=-\varphi^{-1}$,
$\epsilon_{n+1}=\Fib_{n+1}+\Fib_n-\varphi\Fib_{n+1}=-\varphi^{-1}\Fib_{n+1}+\Fib_n=-\varphi^{-1}\epsilon_n$.
So $\epsilon_n=(-\varphi^{-1})^n$, and $\nint{\varphi\Fib_n}\le|\varphi\Fib_n-\Fib_{n+1}|=\varphi^{-n}$. For $n\ge2$,
$\varphi^{-n}\le\varphi^{-2}=2-\varphi<\frac12$, so $\Fib_{n+1}$ is the integer nearest to $\varphi\Fib_n$ and equality
holds; for $n=1$, $\nint{\varphi}=2-\varphi=\varphi^{-2}$.
(ii) Here $u\tau_n=\Fib_n\in\Z$, so $\nint{u\tau_n}=0$, and $\nint{\varphi u\tau_n}\le\varphi^{-n}$; apply
Theorem~\ref{thm:T4}(b).
(iii) Since $\theta\mapsto\cos2\pi\theta$ is even and $1$-periodic:
$\cos(2\pi\Fib_n)=1$; $\cos(2\pi\varphi\Fib_n)=\cos(2\pi(\Fib_{n+1}-\epsilon_n))=\cos(2\pi\varphi^{-n})$;
$\cos(2\pi\varphi^2\Fib_n)=\cos(2\pi(\varphi\Fib_n+\Fib_n))=\cos(2\pi\varphi\Fib_n)$; and
$\cos(2\pi\varphi^3\Fib_n)=\cos(2\pi(2\varphi\Fib_n+\Fib_n))=\cos(4\pi\epsilon_n)=\cos(4\pi\varphi^{-n})$.
Summing with multiplicities $1,2,2,1$ gives the formula. For $n\ge1$, $\tau_n\ne0$, so $C_S(\tau_n)<6$ by
Proposition~\ref{prop:closed}(ii). Finally, with $\vartheta=\varphi^{-n}$ and $1-\cos y\le y^2/2$,
$6-C_S(\tau_n)=4(1-\cos2\pi\vartheta)+(1-\cos4\pi\vartheta)\le8\pi^2\vartheta^2+8\pi^2\vartheta^2$.
\end{proof}

\begin{lemma}[density]\label{lem:dense}
For $u>0$ the set $\{(ut,\varphi ut)\bmod\Z^2:\ t\in\R\}$ is dense in $\T^2$.
\end{lemma}

\begin{proof}
Let $(\xi_1,\xi_2)\in\R^2$ and $\eta>0$. For $t=(\xi_1+k)/u$, $k\in\Z$, the point is $(\xi_1,\varphi\xi_1+k\varphi)\bmod\Z^2$,
so it suffices to show that $\{k\varphi\bmod1:k\in\Z\}$ is dense in $\T$. Choose an integer $N>1/\eta$. The fractional
parts of $k\varphi$, $0\le k\le N$, are $N+1$ distinct numbers in $[0,1)$ (as $\varphi\notin\Q$), so two of them, for
$k_1<k_2$, differ by less than $1/N$; then $m=k_2-k_1\ge1$ satisfies $0<\delta:=\nint{m\varphi}<\eta$. The points
$jm\varphi\bmod1$, $j\in\Z$, are the points $\pm j\delta\bmod1$ (with a fixed sign), so every closed arc of length $\delta$
contains one of them; in particular some $j$ satisfies $\nint{jm\varphi-(\xi_2-\varphi\xi_1)}\le\delta/2<\eta$.
\end{proof}

This is the special case of Kronecker's theorem \cite{HardyWright2008} needed below.

\begin{remark}[torus representation]\label{rem:torus}
With $\Gamma(\theta_1,\theta_2)=\cos2\pi\theta_1+2\cos2\pi\theta_2+2\cos2\pi(\theta_1+\theta_2)+\cos2\pi(\theta_1+2\theta_2)
=8\cos^2(\pi\theta_2)\cos^2(\pi(\theta_1+\theta_2))-2$, a continuous function on $\T^2$, one has
$C_S(t)=\Gamma(ut,\varphi ut)$. The exponents of $C_S$ lie in the $\Z$-module $u\Z+\varphi u\Z=u\,\Z[\varphi]$ of rank $2$.
By Lemma~\ref{lem:dense} the
line $t\mapsto(ut,\varphi ut)$ is dense in $\T^2$, so $C_S$ is quasiperiodic in the classical sense: the
restriction of a continuous function on a two-dimensional torus to a dense line. In particular
$\sup_{t\ge R}C_S(t)=\max \Gamma=6$ for every $R$, by Proposition~\ref{prop:fibAP}(iii).
\end{remark}

The bound $12\pi\varphi^{-n}$ in Proposition~\ref{prop:fibAP}(ii) is not optimal. The torus representation gives
the exact value of the left-hand side.

\begin{proposition}[sharp constant]\label{prop:sharp}
For $n\ge1$ put $\tau_n=\Fib_n/u$, $E_n=\sup_{t\in\R}|C_S(t+\tau_n)-C_S(t)|$ and $\delta_n=\varphi^{-n}$. For
$-1\le\gamma\le1$ let $\mu(\gamma)=\max_{X,Y\in\R}|\sin X+\sin Y+\gamma\sin(X+Y)|$.
\begin{enumerate}[label=\textup{(\roman*)}]
\item $\mu(-\gamma)=\mu(\gamma)$, and for $0\le\gamma\le1$ (cf.\ \cite[(19)--(21)]{Goldobin2007} and Remark~\ref{rem:goldobin})
  \[
    \mu(\gamma)=\max_{0\le s\le1}2\sqrt{1-s^2}\,(1+\gamma s);
  \]
  for $\gamma>0$ the maximum is attained at $s^*(\gamma)=\bigl(\sqrt{1+8\gamma^2}-1\bigr)/(4\gamma)$, and $\mu(0)=2$.
  The function $\mu$ is non-decreasing on $[0,1]$, and $\mu(1)=\tfrac{3\sqrt3}{2}$, attained at $X=Y=\pi/3$.
\item $E_n=4\sin(\pi\delta_n)\,\mu(\cos\pi\delta_n)$ for all $n\ge1$. In particular $E_n\le6\sqrt3\,\pi\varphi^{-n}$.
\item For $n\ge2$,
  \[
    E_n=6\sqrt3\,\pi\varphi^{-n}-2\sqrt3\,\pi^3\varphi^{-3n}+R_n,\qquad
    \tfrac{\sqrt3}{6}\,\pi^5\varphi^{-5n}\le R_n\le\Bigl(\tfrac38+\tfrac{3\sqrt3}{10}\Bigr)\pi^5\varphi^{-5n}.
  \]
\item $R_n=\frac{16\sqrt3}{45}\,\pi^5\varphi^{-5n}+O(\varphi^{-7n})$.
\end{enumerate}
In particular the leading constant $6\sqrt3\,\pi$ and the second-order constant $2\sqrt3\,\pi^3$ are both optimal, and
$6\sqrt3\,\pi=4\pi\mu(1)$ is the norm $\sup_{\T^2}|\partial_{\theta_2}\Gamma|$ of the derivative of $\Gamma$ in the
direction of the shift.
\end{proposition}

\begin{proof}
(i) The substitution $(X,Y)\mapsto(X+\pi,Y+\pi)$ changes the signs of $\sin X$ and $\sin Y$ and leaves $\sin(X+Y)$
unchanged, so $\mu(-\gamma)=\mu(\gamma)$. Let $0\le\gamma\le1$ and put $\Sigma=X+Y$, $D=(X-Y)/2$; as $(X,Y)$ ranges over
$\R^2$ so does $(\Sigma,D)$. Since $\sin X+\sin Y=2\sin(\Sigma/2)\cos D$ and $\sin\Sigma=2\sin(\Sigma/2)\cos(\Sigma/2)$,
\[
  \sin X+\sin Y+\gamma\sin(X+Y)=2\sin\tfrac{\Sigma}{2}\Bigl(\cos D+\gamma\cos\tfrac{\Sigma}{2}\Bigr).
\]
For fixed $\Sigma$ the maximum of the modulus over $D$ is $2|\sin(\Sigma/2)|\,(1+\gamma|\cos(\Sigma/2)|)$ (take $\cos D=\pm1$
with the sign of $\cos(\Sigma/2)$), and $s=|\cos(\Sigma/2)|$ ranges over $[0,1]$ with $|\sin(\Sigma/2)|=\sqrt{1-s^2}$. This
proves the formula; each function $\gamma\mapsto2\sqrt{1-s^2}(1+\gamma s)$ is non-decreasing, hence so is $\mu$. For
$\gamma>0$ and $0\le s<1$, $\frac{d}{ds}\bigl[\sqrt{1-s^2}(1+\gamma s)\bigr]=(\gamma-s-2\gamma s^2)/\sqrt{1-s^2}$, which is
positive for $s<s^*(\gamma)$ and negative for $s>s^*(\gamma)$, where $s^*(\gamma)\in(0,1)$ is the positive root of
$2\gamma s^2+s-\gamma=0$. For $\gamma=1$, $s^*=\frac12$ and $\mu(1)=2\cdot\frac{\sqrt3}2\cdot\frac32=\frac{3\sqrt3}2$; $s=\frac12$
and $D=0$ correspond to $\Sigma=2\pi/3$, $X=Y=\pi/3$.

(ii) Let $\epsilon_n=\Fib_{n+1}-\varphi\Fib_n=(-\varphi^{-1})^n$ (Proposition~\ref{prop:fibAP}(i)) and $\delta=-\epsilon_n$, so
$|\delta|=\delta_n$. With $\Gamma$ as in Remark~\ref{rem:torus} and $(\theta_1,\theta_2)=(ut,\varphi ut)$ we have
$u(t+\tau_n)=\theta_1+\Fib_n$ and $\varphi u(t+\tau_n)=\theta_2+\Fib_{n+1}+\delta$, so by $1$-periodicity of $\Gamma$ in each
variable, $C_S(t+\tau_n)-C_S(t)=\Gamma(\theta_1,\theta_2+\delta)-\Gamma(\theta_1,\theta_2)$. The variable $\theta_2$ enters the
four terms of $\Gamma$ with coefficients $0,1,1,2$. Using $\cos(\alpha+2k\pi\delta)-\cos\alpha=-2\sin(k\pi\delta)\sin(\alpha+k\pi\delta)$
and $\sin2\pi\delta=2\sin\pi\delta\cos\pi\delta$, and putting $X=2\pi\theta_2+\pi\delta$, $Y=2\pi(\theta_1+\theta_2)+\pi\delta$
(so that $2\pi(\theta_1+2\theta_2)+2\pi\delta=X+Y$), we obtain the identity
\[
  \Gamma(\theta_1,\theta_2+\delta)-\Gamma(\theta_1,\theta_2)=-4\sin(\pi\delta)\bigl(\sin X+\sin Y+\cos(\pi\delta)\sin(X+Y)\bigr).
\]
The map $(\theta_1,\theta_2)\mapsto(X,Y)$ is affine with linear part $2\pi\bigl(\begin{smallmatrix}0&1\\1&1\end{smallmatrix}\bigr)$,
and the integer matrix $\bigl(\begin{smallmatrix}0&1\\1&1\end{smallmatrix}\bigr)$ has determinant $-1$, so the map induces
a homeomorphism $\T^2\to(\R/2\pi\Z)^2$. By Lemma~\ref{lem:dense} the image of $t\mapsto(\theta_1,\theta_2)$ is dense,
hence so is the image of $t\mapsto(X,Y)$. The function $(X,Y)\mapsto|\sin X+\sin Y+\gamma\sin(X+Y)|$ is continuous on
the compact torus $(\R/2\pi\Z)^2$, so its supremum over a dense subset is its maximum $\mu(\gamma)$. Since
$0<\delta_n<1$, $|\sin\pi\delta|=\sin\pi\delta_n$, and $\cos\pi\delta=\cos\pi\delta_n$; this proves the formula. By (i),
$\mu(\cos\pi\delta_n)\le\mu(1)=\frac{3\sqrt3}2$, and $\sin y\le y$ gives $E_n\le6\sqrt3\,\pi\delta_n$.

(iii) Let $n\ge2$, $y=\pi\delta_n\le\pi\varphi^{-2}<1.2$ and $\gamma=\cos y\in(0,1)$. We first show
\begin{equation}\label{eq:mu-two-sided}
  \tfrac{3\sqrt3}{2}-\tfrac{\sqrt3}{2}(1-\gamma)\;\le\;\mu(\gamma)\;\le\;\tfrac{3\sqrt3}{2}-\tfrac{\sqrt3}{2}(1-\gamma)+\tfrac38(1-\gamma)^2 .
\end{equation}
The lower bound is the value at $s=\frac12$ in (i). For the upper bound write
$2\sqrt{1-s^2}(1+\gamma s)=f(s)-(1-\gamma)g(s)$ with $f(s)=2\sqrt{1-s^2}(1+s)$ and $g(s)=2s\sqrt{1-s^2}$. On $[0,1]$ we have
$\frac{3\sqrt3}2-2(s-\frac12)^2\ge\frac{3\sqrt3}2-\frac12>0$ and
\[
  \Bigl(\tfrac{3\sqrt3}{2}-2\bigl(s-\tfrac12\bigr)^2\Bigr)^2-f(s)^2=\bigl(s-\tfrac12\bigr)^2\bigl(8s^2+8s+12-6\sqrt3\bigr)\ge0,
\]
so $f(s)\le\frac{3\sqrt3}2-2(s-\frac12)^2$. Next, $g''(s)=2s(2s^2-3)(1-s^2)^{-3/2}\le0$ on $[0,1)$, so $g$ is concave on
$[0,1]$; as $g(0)=g(1)=0$ and $g(\frac12)=\frac{\sqrt3}2$, concavity gives $g(s)\ge\frac{\sqrt3}2-\sqrt3\,|s-\frac12|$. With
$d=|s-\frac12|$,
\[
  2\sqrt{1-s^2}(1+\gamma s)\le\tfrac{3\sqrt3}{2}-\tfrac{\sqrt3}{2}(1-\gamma)+\bigl(\sqrt3(1-\gamma)d-2d^2\bigr)
  \le\tfrac{3\sqrt3}{2}-\tfrac{\sqrt3}{2}(1-\gamma)+\tfrac38(1-\gamma)^2,
\]
which is \eqref{eq:mu-two-sided}. Now use $y-\frac{y^3}6\le\sin y\le y-\frac{y^3}6+\frac{y^5}{120}$ and
$\frac{y^2}2-\frac{y^4}{24}\le1-\cos y\le\frac{y^2}2$ (valid for $y\ge0$). By (ii) and \eqref{eq:mu-two-sided}, since all
factors below are positive for $0<y<1.2$,
\[
  E_n\ge4\Bigl(y-\tfrac{y^3}{6}\Bigr)\Bigl(\tfrac{3\sqrt3}{2}-\tfrac{\sqrt3}{4}y^2\Bigr)=6\sqrt3\,y-2\sqrt3\,y^3+\tfrac{\sqrt3}{6}y^5,
\]
\[
  E_n\le4\Bigl(y-\tfrac{y^3}{6}+\tfrac{y^5}{120}\Bigr)\Bigl(\tfrac{3\sqrt3}{2}-\tfrac{\sqrt3}{4}y^2+\Bigl(\tfrac{\sqrt3}{48}+\tfrac{3}{32}\Bigr)y^4\Bigr)
  =6\sqrt3\,y-2\sqrt3\,y^3+q(y)\,y^5,
\]
where $q(y)=\frac38+\frac{3\sqrt3}{10}-\bigl(\frac1{16}+\frac{\sqrt3}{45}\bigr)y^2+\bigl(\frac1{320}+\frac{\sqrt3}{1440}\bigr)y^4\le\frac38+\frac{3\sqrt3}{10}$
for $0<y\le4$. With $y=\pi\varphi^{-n}$ this is (iii). Optimality of the two constants follows from (iii). Finally,
$\partial_{\theta_2}\Gamma=-4\pi\bigl(\sin2\pi\theta_2+\sin2\pi(\theta_1+\theta_2)+\sin2\pi(\theta_1+2\theta_2)\bigr)$ has supremum norm
$4\pi\mu(1)=6\sqrt3\,\pi$ by (i).

(iv) Near $\gamma=1$, $s^*(\gamma)$ is a simple root of $2\gamma s^2+s-\gamma=0$, hence real-analytic, and
$\mu(\gamma)=\Psi(s^*(\gamma),\gamma)$ with $\Psi(s,\gamma)=2\sqrt{1-s^2}(1+\gamma s)$ and $\partial_s\Psi(s^*,\gamma)=0$. Hence
$\mu'=2s^*\sqrt{1-s^{*2}}$ and $\mu''=2s^{*\prime}(1-2s^{*2})/\sqrt{1-s^{*2}}$. At $\gamma=1$: $s^*=\frac12$ and, differentiating
$2\gamma s^{*2}+s^*-\gamma=0$, $s^{*\prime}(1)=(1-2s^{*2})/(4\gamma s^*+1)=\frac16$; so $\mu'(1)=\frac{\sqrt3}2$ and
$\mu''(1)=\frac{\sqrt3}9$. With $1-\cos y=\frac{y^2}2-\frac{y^4}{24}+O(y^6)$,
\[
  \mu(\cos y)=\tfrac{3\sqrt3}{2}-\tfrac{\sqrt3}{4}y^2+\tfrac{5\sqrt3}{144}y^4+O(y^6),
\]
and multiplying by $4\sin y=4y-\frac23y^3+\frac1{30}y^5+O(y^7)$ gives
$E_n=6\sqrt3\,y-2\sqrt3\,y^3+\frac{16\sqrt3}{45}y^5+O(y^7)$ with $y=\pi\varphi^{-n}$.
\end{proof}

\begin{remark}\label{rem:goldobin}
Writing $s=\cos\theta$, $2\sqrt{1-s^2}\,(1+\gamma s)=2\bigl(\sin\theta+\tfrac{\gamma}{2}\sin2\theta\bigr)$,
so $\mu(\gamma)=2\mathcal M(\gamma/2)$ with $\mathcal M(g)=\max_\theta(\sin\theta+g\sin2\theta)$.
This maximum was computed in closed form in \cite[(19)--(21)]{Goldobin2007} (their $\gamma_{c+}(g)$):
the maximiser satisfies $\cos\theta_0=(\sqrt{1+32g^2}-1)/(8g)$, which is $s^*(\gamma)$ for
$g=\gamma/2$, and for $0<\gamma\le1$, $r=\sqrt{1+8\gamma^2}$, one gets
$\mu(\gamma)=(r+3)^{3/2}(r-1)^{1/2}/(8\gamma)$; e.g.\ $\mu(1)=\sqrt{432}/8=3\sqrt3/2$.
\end{remark}

\begin{remark}
The upper bound $E_n\le6\sqrt3\,\pi\varphi^{-n}$ is the mean value estimate $|\delta|\sup|\partial_{\theta_2}\Gamma|$;
Proposition~\ref{prop:sharp} shows that it is asymptotically attained, and also identifies the next term. It improves
the bound $12\pi\varphi^{-n}$ of Proposition~\ref{prop:fibAP}(ii). By (iv), $R_n\varphi^{5n}/\pi^5\to16\sqrt3/45=0.61584\ldots$
(Section~\ref{sec:numerics}), inside the interval $[0.289,\,0.895]$ given by (iii).
\end{remark}

For $\varepsilon>0$ let $L(\varepsilon)$ be the infimum of the numbers $L>0$ such that every interval of length $L$ contains
an element of $E(C_S;\varepsilon)$ (Definition~\ref{def:ap}).

\begin{proposition}[order of $L(\varepsilon)$]\label{prop:Leps}
For $0<\varepsilon\le\frac12$,
\[
  \frac{\sqrt2}{3\varphi^2u}\cdot\frac1\varepsilon\;\le\;L(\varepsilon)\;\le\;\frac1u\Bigl(\frac{6\pi\varphi^3}{\varepsilon}+1\Bigr).
\]
\end{proposition}

\begin{proof}
\emph{Upper bound.} Let $\varepsilon\le1$. Put $\eta=\varepsilon/(12\pi)$ and let $m$ be the least integer with
$\varphi^{-m+1}<2\eta$; then $m\ge2$ and $\varphi^{-m+2}\ge2\eta$, i.e.\ $\varphi^m\le\varphi^2/(2\eta)$. Consider the
lattice $\Lambda=\{(k,k\varphi-j):k,j\in\Z\}\subset\R^2$, the image of $\Z^2$ under $(k,j)\mapsto(k,k\varphi-j)$. By
Cassini's identity (Lemma~\ref{lem:cassini}), $(\Fib_m,\Fib_{m+1})$ and $(\Fib_{m+1},\Fib_{m+2})$ form a basis of $\Z^2$, so
$v_1=(\Fib_m,-\epsilon_m)$ and $v_2=(\Fib_{m+1},-\epsilon_{m+1})$, with $\epsilon_m=\Fib_{m+1}-\varphi\Fib_m=(-\varphi^{-1})^m$, form a
basis of $\Lambda$. Hence every translate of the half-open parallelogram $P=\{sv_1+tv_2:0\le s,t<1\}$ contains exactly
one point of $\Lambda$. The first coordinates of $P$ fill $[0,\Fib_{m+2})$, and, since $\epsilon_m$ and $\epsilon_{m+1}$ have
opposite signs, the second coordinates fill an interval of length $|\epsilon_m|+|\epsilon_{m+1}|=\varphi^{-m}+\varphi^{-m-1}=\varphi^{-m+1}<2\eta$.
Hence for every integer $k_0$ the box $[k_0,k_0+\Fib_{m+2})\times(-\eta,\eta)$ contains a translate of $P$, and therefore a
point of $\Lambda$: there is an integer $k\in[k_0,k_0+\Fib_{m+2})$ with $\nint{k\varphi}<\eta$. For such $k$, $\tau=k/u$ has
$\nint{u\tau}=0$ and $\nint{\varphi u\tau}<\eta$, so $\sup_t|C_S(t+\tau)-C_S(t)|\le12\pi\eta=\varepsilon$ by Theorem~\ref{thm:T4}(b),
i.e.\ $\tau\in E(C_S;\varepsilon)$. Every interval of length $(\Fib_{m+2}+1)/u$ contains $\Fib_{m+2}$ consecutive points $k/u$, and
$\Fib_{m+2}\le\varphi^{m+1}$ (by induction: $\Fib_1\le\varphi^0$, $\Fib_2\le\varphi$, $\varphi^{k+1}=\varphi^k+\varphi^{k-1}$), and $\varphi^{m+1}\le\varphi^3/(2\eta)=6\pi\varphi^3/\varepsilon$. This gives the upper bound.

\emph{Lower bound.} For $\omega\in\R^2$ put $\mathcal F(\omega)=\sup_{\theta\in\R^2}|\Gamma(\theta+\omega)-\Gamma(\theta)|$. By
Lemma~\ref{lem:dense} (and continuity of $\theta\mapsto|\Gamma(\theta+\omega)-\Gamma(\theta)|$),
$\sup_t|C_S(t+\tau)-C_S(t)|=\mathcal F(u\tau,\varphi u\tau)$. Write $\Gamma(\theta)=\sum_{\kappa\in\mathcal K_0}r_\kappa\cos(2\pi\kappa\cdot\theta)$
with $\mathcal K_0=\{(1,0),(0,1),(1,1),(1,2)\}$ and $(r_\kappa)=(1,2,2,1)$. Then
\[
  \Gamma(\theta+\omega)-\Gamma(\theta)=-2\sum_{\kappa\in\mathcal K_0}r_\kappa\sin(\pi\kappa\cdot\omega)\,\sin\bigl(2\pi\kappa\cdot\theta+\pi\kappa\cdot\omega\bigr).
\]
No two elements of $\mathcal K_0\cup(-\mathcal K_0)$ coincide, so these four functions of $\theta$ are orthogonal in
$L^2(\T^2)$, each of norm $2^{-1/2}$. Since $\sup|\psi|\ge\|\psi\|_{L^2(\T^2)}$ and $|\sin\pi x|\ge2\nint{x}$ (concavity of
$\sin$ on $[0,\pi/2]$),
\[
  \mathcal F(\omega)^2\ge2\sum_{\kappa\in\mathcal K_0}r_\kappa^2\sin^2(\pi\kappa\cdot\omega)\ge2\sin^2(\pi\omega_1)+8\sin^2(\pi\omega_2)
  \ge8\nint{\omega_1}^2+32\nint{\omega_2}^2\ge8\nint{\omega}^2,
\]
where $\nint{\omega}=\max(\nint{\omega_1},\nint{\omega_2})$. Thus $\mathcal F(\omega)\ge c\nint\omega$ with $c=2\sqrt2$.

Let $\tau\in E(C_S;\varepsilon)$, $\tau>0$. Then $\nint{u\tau},\nint{\varphi u\tau}\le\varepsilon/c$. Write $u\tau=k+s$ with $k\in\Z$ and
$|s|\le\varepsilon/c\le1/(4\sqrt2)<1$; since $u\tau>0$, $k\ge0$. If $k=0$, then $u\tau\le\varepsilon/c$. If $k\ge1$, then
$\nint{k\varphi}\le\nint{\varphi u\tau}+\varphi|s|\le\varphi^2\varepsilon/c$, while $\nint{k\varphi}>1/(3k)$: if $j$ is the integer nearest to
$k\varphi$, then $(j-k\varphi)(j-k\bar\varphi)=j^2-jk-k^2$ is a non-zero integer and $|j-k\bar\varphi|=|j-k\varphi+k\sqrt5|\le\frac12+\sqrt5k<3k$.
Hence $k>c/(3\varphi^2\varepsilon)$ and $u\tau\ge k-\varepsilon/c>c/(3\varphi^2\varepsilon)-\varepsilon/c$. Since
$\varepsilon\le\frac12<\sqrt{2/3}/\varphi=c/(2\sqrt3\,\varphi)$, we have $2\varepsilon/c\le c/(6\varphi^2\varepsilon)$, i.e.\
$c/(3\varphi^2\varepsilon)-2\varepsilon/c\ge c/(6\varphi^2\varepsilon)$. So the interval
$\bigl(\varepsilon/(cu),\,(c/(3\varphi^2\varepsilon)-\varepsilon/c)/u\bigr)$ contains no element of $E(C_S;\varepsilon)$ and has length at
least $c/(6\varphi^2u\varepsilon)=\sqrt2/(3\varphi^2u\varepsilon)$; hence $L(\varepsilon)\ge\sqrt2/(3\varphi^2u\varepsilon)$.
\end{proof}

\begin{remark}
The order $1/\varepsilon$ reflects the fact that $\varphi$ is badly approximable. For the return times of $\{n\theta\}$ to an interval (Slater's three-gap theorem)
see \cite{Slater1967}, and for quasi-uniformity of Kronecker sequences with badly
approximable $\theta$ see Goda \cite{Goda2024}. Upper bounds for the inclusion length via simultaneous Diophantine approximation
are due to Naito \cite{Naito1996}; Anikushin \cite[Theorems~3 and~4]{Anikushin2017} gives a lower bound for the growth exponent of the
inclusion length of quasiperiodic functions and notes that both agree for badly approximable exponents. The
contribution of Proposition~\ref{prop:Leps} is the two-sided inequality for this particular $C_S$, with explicit
constants, on the whole range $0<\varepsilon\le\frac12$.
\end{remark}

The method (lattice points in a thin strip, and the bounded partial quotients of $\varphi$) is standard; we have not
found the statement for $C_S$ in the literature.

\section{Numerical examples}\label{sec:numerics}

The tabulated numerical values were computed in $60$-digit floating-point arithmetic and rounded.

\begin{example}[$u=1$, $F=\varphi^4$]\label{ex:unit}
Here $a=\varphi$, $b=\varphi^2$, $a+b=\varphi^3\approx4.236\,068<F=\varphi^4\approx6.854\,102$. By
Lemma~\ref{lem:phi}(ii),
\[
  S=(\varphi^4,\ \varphi^4-\varphi,\ \varphi^4-\varphi^2,\ \varphi^4-\varphi^3)=(3\varphi+2,\ 2\varphi+2,\ 2\varphi+1,\ \varphi+1)
   =(\varphi^4,\ 2\varphi^2,\ \varphi^3,\ \varphi^2),
\]
numerically $(6.854\,102,\ 5.236\,068,\ 4.236\,068,\ 2.618\,034)$. Table~\ref{tab:unit} lists the six
differences exactly in $\Z+\Z\varphi$.
\end{example}

\begin{table}[ht]
\centering
\caption{Differences of $S$ for $u=1$, $F=\varphi^4$ (Example~\ref{ex:unit}).}\label{tab:unit}
\begin{tabular}{@{}lccc@{}}
\toprule
pair & exact value in $\Z+\Z\varphi$ & power of $\varphi$ & decimal\\
\midrule
$s_1-s_2$ & $\varphi$        & $\varphi^1$ & $1.618\,034$\\
$s_1-s_3$ & $\varphi+1$      & $\varphi^2$ & $2.618\,034$\\
$s_1-s_4$ & $2\varphi+1$     & $\varphi^3$ & $4.236\,068$\\
$s_2-s_3$ & $1$              & $\varphi^0$ & $1.000\,000$\\
$s_2-s_4$ & $\varphi+1$      & $\varphi^2$ & $2.618\,034$\\
$s_3-s_4$ & $\varphi$        & $\varphi^1$ & $1.618\,034$\\
\bottomrule
\end{tabular}
\end{table}

\begin{example}[the Planck scale, $u=\nu_{203}$, $F=\nu_{199}$]\label{ex:planck2}
Let $K=\Kp$ as in Example~\ref{ex:planck}, $u=\nu_{203}$ and $F=\nu_{199}>\nu_{200}=a+b$. By
Lemma~\ref{lem:ladder}, $\nu_{199}-\nu_{200}=\nu_{201}$, $\nu_{199}-\nu_{201}=\nu_{200}$, and
$\nu_{199}-\nu_{202}=(\varphi^4-\varphi)u=(2\varphi+2)u=2\varphi^2u=2\nu_{201}$. Hence
\[
  S=(\nu_{199},\ 2\nu_{201},\ \nu_{200},\ \nu_{201}),
\]
and by Corollary~\ref{cor:ladder} the differences are $\nu_{202},\nu_{201},\nu_{200},\nu_{203},\nu_{201},\nu_{202}$.
Table~\ref{tab:planck} gives the values. They are computed from the central CODATA values; because of the
relative uncertainty $1.1\times10^{-5}$ of $\Kp$, only the first five significant digits are meaningful as
consequences of the measured constants.
\end{example}

\begin{table}[ht]
\centering
\caption{The golden ladder at the Planck scale, $\nu_n=\Kp\varphi^{-n}$ (Example~\ref{ex:planck2}).}\label{tab:planck}
\begin{tabular}{@{}lll@{}}
\toprule
quantity & expression & value\\
\midrule
$\nu_{203}=u$            & $\Kp\varphi^{-203}$ & $2.784\,39$\\
$\nu_{202}=a$            & $\Kp\varphi^{-202}$ & $4.505\,24$\\
$\nu_{201}=b$            & $\Kp\varphi^{-201}$ & $7.289\,63$\\
$\nu_{200}=a+b$          & $\Kp\varphi^{-200}$ & $11.794\,86$\\
$\nu_{199}=F$            & $\Kp\varphi^{-199}$ & $19.084\,49$\\
\midrule
$S$ & $(\nu_{199},2\nu_{201},\nu_{200},\nu_{201})$ & $(19.084\,49,\ 14.579\,25,\ 11.794\,86,\ 7.289\,63)$\\
$\Dl(S)$ & $\{\nu_{203},\nu_{202},\nu_{201},\nu_{200}\}$ & $\{2.784\,39,\ 4.505\,24,\ 7.289\,63,\ 11.794\,86\}$\\
\bottomrule
\end{tabular}
\end{table}

\begin{table}[ht]
\centering
\caption{Closure residues (in units of $u$) for several bases $x$ (Theorem~\ref{thm:T3},
Proposition~\ref{prop:fibres}). Fibonacci quotients are given exactly.}\label{tab:residues}
\begin{tabular}{@{}lccc@{}}
\toprule
$x$ & decimal & $\rho_-(x)/u=x^2-x-1$ & $\rho_+(x)/u=x(x^2-x-1)$\\
\midrule
$\sqrt2$        & $1.414\,214$ & $1-\sqrt2\approx-0.414\,214$ & $\sqrt2-2\approx-0.585\,786$\\
$3/2$           & $1.500\,000$ & $-1/4$     & $-3/8$\\
$8/5$           & $1.600\,000$ & $-1/25$    & $-8/125$\\
$21/13$         & $1.615\,385$ & $-1/169$   & $-21/2197$\\
$\varphi$       & $1.618\,034$ & $0$        & $0$\\
$34/21$         & $1.619\,048$ & $1/441$    & $34/9261$\\
$13/8$          & $1.625\,000$ & $1/64$     & $13/512$\\
$5/3$           & $1.666\,667$ & $1/9$      & $5/27$\\
$\sqrt3$        & $1.732\,051$ & $2-\sqrt3\approx0.267\,949$ & $2\sqrt3-3\approx0.464\,102$\\
$2$             & $2.000\,000$ & $1$        & $2$\\
\bottomrule
\end{tabular}
\end{table}

\begin{table}[ht]
\centering
\caption{Explicit almost periods $\tau_n=\Fib_n$ of $C_S$ for $u=1$ (Proposition~\ref{prop:fibAP}).}\label{tab:AP}
\begin{tabular}{@{}rrcccc@{}}
\toprule
$n$ & $\Fib_n$ & $C_S(\Fib_n)$ & $6-C_S(\Fib_n)$ & $16\pi^2\varphi^{-2n}$ & $12\pi\varphi^{-n}$\\
\midrule
5  & 5    & $4.798\,867$   & $1.201\,1$                & $1.283\,9$                & $3.399\,3$\\
8  & 21   & $5.928\,715$   & $7.128\,5\times10^{-2}$   & $7.155\,1\times10^{-2}$   & $0.802\,47$\\
10 & 55   & $5.989\,566$   & $1.043\,4\times10^{-2}$   & $1.043\,9\times10^{-2}$   & $0.306\,52$\\
12 & 144  & $5.998\,477$   & $1.522\,9\times10^{-3}$   & $1.523\,1\times10^{-3}$   & $0.117\,08$\\
15 & 610  & $5.999\,915$   & $8.487\,7\times10^{-5}$   & $8.487\,7\times10^{-5}$   & $2.763\,9\times10^{-2}$\\
20 & 6765 & $5.999\,999\,3$& $6.901\,0\times10^{-7}$   & $6.901\,0\times10^{-7}$   & $2.492\,2\times10^{-3}$\\
\bottomrule
\end{tabular}
\end{table}

\begin{remark}
For $u=1$ and $n=8,12,16$ a sampling of $|C_S(t+\Fib_n)-C_S(t)|$ on a grid of $4000$ points in $[-50,50]$ gave maxima
$0.694$, $0.101$ and $0.0148$, against the bounds $12\pi\varphi^{-n}=0.802$, $0.117$ and $0.0171$ of
Proposition~\ref{prop:fibAP}(ii). This is only a numerical illustration and plays no role in the proofs.
\end{remark}

\begin{remark}[check of Proposition~\ref{prop:sharp}]
For $u=1$ we compared, for $n\in\{1,2,5,10,15\}$, the exact value $4\sin(\pi\varphi^{-n})\mu(\cos\pi\varphi^{-n})$ of
Proposition~\ref{prop:sharp}(ii) with the maximum of $|C_S(t+\Fib_n)-C_S(t)|$ over $2\cdot10^6$ random $t\in[0,2\cdot10^4]$;
the sampled maxima lie below the exact values, with relative differences between $2\cdot10^{-8}$ and $10^{-6}$ for this particular sample (other random samples give relative differences up to about $7\cdot10^{-6}$),
as expected from finite sampling (for $n=10$: $0.2653938$ sampled, $0.2653939$ exact). Independently, the closed
form of $\mu$ in (i) agrees with a direct maximisation over the torus to $10^{-15}$ for $21$ values of
$\gamma\in[0,1]$. The difference between the exact value and $6\sqrt3\,\pi\varphi^{-n}$ is, relatively, $1.4\times10^{-8}$
for $n=20$; and $R_n\varphi^{5n}/\pi^5$, computed in $50$-digit arithmetic, equals $0.5632$, $0.6069$, $0.6158$ and
$0.61584$ for $n=2,4,10,20$, within the bounds of Proposition~\ref{prop:sharp}(iii). This check plays no role in
the proofs.
\end{remark}

\section{Relation to the literature: what is classical and what is new}\label{sec:literature}

\subsection{Sidon sets and Golomb rulers}
Sidon sets go back to Sidon's work on trigonometric series \cite{Sidon1932}; the basic counting bound is due
to Erd\H{o}s and Tur\'an \cite{ErdosTuran1941}, and O'Bryant's survey \cite{OBryant2004} is a comprehensive guide.
Integer Sidon sets are known as Golomb rulers \cite{Golomb1972,Drakakis2009}. The golden quadruple is
\emph{not} a Sidon set (Remark~\ref{rem:structure}(c)), so it is not a ruler in this sense. There is,
however, a formal parallel: a perfect Golomb ruler (Definition~\ref{def:classical}) has at most four marks,
e.g.\ $\{0,1,4,6\}$ with differences $\{1,\dots,6\}$. After a translation to $\min A=0$, the marks of a perfect
ruler with $n$ marks are exactly the vertex labels of a graceful labelling of the complete graph $K_n$ (an injective
map $V(K_n)\to\{0,1,\dots,\binom n2\}$ whose induced edge labels $|f(v)-f(w)|$ are pairwise distinct), so this is
the classical fact that $K_n$ is graceful if and only if $n\le4$, in the graph formulation of Golomb \cite{Golomb1972}.
For the equivalent problem of zero-redundancy linear arrays, Moffet \cite{Moffet1968} credits the
construction to Arsac (1955) and the proof that these are the only such arrays to Bracewell (1962); see also
\cite[\S2.4]{Gallian2025}.
Theorems~\ref{thm:T2} and~2$'$ exhibit a four-point barrier when the prescribed set of differences is the
geometric progression $\Dset$, or the orbit $\{c\varphi^n\}$, instead of $\{1,\dots,L\}$. The two statements are
logically independent.

\subsection{Difference sets and the size of $A-A$}
Combinatorial difference sets, beginning with Singer's construction \cite{Singer1938}, are treated in
\cite{Baumert1971,BJL1999}. For finite $A\subset\R$ the bounds $2|A|-1\le|A-A|\le|A|^2-|A|+1$, with equality
exactly for arithmetic progressions and for Sidon sets respectively, are standard \cite{Nathanson1996,TaoVu2006}.
The golden quadruple has $|S-S|=9$ and is a translate of the two-dimensional cube $\{0,a\}+\{0,b\}$; such cubes
are often called Hilbert cubes \cite{Hilbert1892,TaoVu2006}. Lemma~\ref{lem:cube} is elementary and classical
in spirit.

\subsection{Fibonacci numbers and the golden ratio}
The identities used here (the recursion of $\varphi^n$, Cassini's identity, $\Fib_{n+1}-\varphi\Fib_n=(-\varphi^{-1})^n$)
are classical \cite{GKP1994,Vorobiev2002,HardyWright2008}. The quality of the almost periods $\Fib_n/u$ in
Proposition~\ref{prop:fibAP} reflects the fact that $\varphi$ is badly approximable, $\nint{\varphi\Fib_n}=\varphi^{-n}$
for $n\ge2$ (Proposition~\ref{prop:fibAP}(i)), and in fact extremal in Hurwitz's theorem \cite{Hurwitz1891,HardyWright2008};
Proposition~\ref{prop:sharp} gives the sharp constant for the specific almost periods $\Fib_n/u$.
A related classical phenomenon is the three-distance theorem of S\'os, \'Swierczkowski and Sur\'anyi
\cite{Sos1958,Swierczkowski1958,Suranyi1958}: for irrational $\alpha\in(0,1)$ and $n\ge1$, the points $\{i\alpha\}$,
$0\le i\le n$, divide the circle $\R/\Z$ into $n+1$ arcs whose lengths take at most three values; when three values occur,
the largest is the sum of the other two. In the formulation of \cite{AlessandriBerthe1998}: if $\alpha=[0;c_1,c_2,\dots]$
has convergents $p_k/q_k$, $\eta_k=(-1)^k(q_k\alpha-p_k)$, and $n=mq_k+q_{k-1}+r$ with $1\le m\le c_{k+1}$ and
$0\le r<q_k$, then the possible arc lengths are $\eta_k$, $\eta_{k-1}-m\eta_k$ and $\eta_{k-1}-(m-1)\eta_k$. For
$\alpha=\varphi^{-1}=[0;1,1,1,\dots]$ all $c_j=1$, hence $m=1$, and $p_k=\Fib_k$, $q_k=\Fib_{k+1}$ for $k\ge0$; since
$\varphi^{-1}=\varphi-1$, Proposition~\ref{prop:fibAP}(i) gives $q_k\alpha-p_k=\varphi\Fib_{k+1}-\Fib_{k+2}=(-1)^k\varphi^{-k-1}$,
i.e.\ $\eta_k=\varphi^{-k-1}$ (and $\eta_{-1}=1$). Thus the possible arc lengths belong to
$\{\varphi^{-k-2},\varphi^{-k-1},\varphi^{-k}\}$.
When all three occur, they are consecutive powers of $\varphi^{-1}$,
with the largest equal to the sum of the other two.
If $r=q_k-1$, only the two smaller lengths occur. Our configuration
has the same flavour: distances that are consecutive powers of $\varphi$, the largest being the sum of two smaller
ones.

\subsection{Quasicrystals}
Since the discovery of Shechtman, Blech, Gratias and Cahn \cite{Shechtman1984}, the golden ratio has been central in
the mathematics of aperiodic order: in de Bruijn's algebraic theory of Penrose tilings \cite{deBruijn1981}, in
Meyer's harmonious sets \cite{Meyer1972,Lagarias1996}, and in the theory of model sets
\cite{Moody1997,BaakeGrimm2013}. The basic one-dimensional example is the Fibonacci chain, generated by the
substitution $L\mapsto LS$, $S\mapsto L$ with tiles of lengths $\varphi$ and $1$ \cite{BaakeGrimm2013}. Since the
second iterate of $L$ is $LSL$, the endpoints of a patch $LSL$ form the set $\{0,\varphi,\varphi^2,\varphi^3\}$, which by
Remark~\ref{rem:structure}(a) is a dilate-translate of the golden quadruple; Theorem~2$'$(ii) shows that it is
the only four-point set all of whose distances are powers of $\varphi$ (the case $c=1$), up to translation and
scaling by $\varphi^n$ (for fixed $c$). The exponents of $C_S$ lie in the rank-two
module $u\Z[\varphi]$ (Remark~\ref{rem:torus}). For the Fibonacci chain with tile lengths $\varphi$ and $1$, both the
return module $\Z[\varphi]$ \cite[\S3]{BaakeGahlerMazac2024} and the Fourier module $\Z[\varphi]/\sqrt5$, which
supports the diffraction measure and is also the Fourier--Bohr spectrum \cite[\S8]{BaakeGahlerMazac2024}, are
modules of this type, namely real multiples of $\Z[\varphi]$; see \cite{BaakeGrimm2013} for the general theory of
model sets.

\subsection{Almost periodic functions}
Almost periodic functions were introduced by Bohr \cite{Bohr1925,Bohr1947}; see also
\cite{Besicovitch1932,Corduneanu1989,Katznelson2004}. That trigonometric polynomials are almost periodic is part of
the foundations of the theory, and the compactness proof of Lemma~\ref{lem:trigAP} is standard. The periodicity
criterion of Proposition~\ref{prop:periodic} is well known; it is a special case of a theorem of Mirotin and
Mirotin on sums of continuous periodic functions \cite[\S2]{MirotinMirotin2009}, and the direct proof given above
uses only Lemma~\ref{lem:square}.

\subsection{Scope: other ratios}
Some of the phenomena above are specific to $\varphi$. Let $\varrho\approx1.324\,718$ be the real root of
$\varrho^3=\varrho+1$ (the plastic number, the smallest Pisot number \cite{Siegel1944,Bertin1992}). Then
$\varrho^{n+3}=\varrho^{n+1}+\varrho^n$ for $n\in\Z$, $\varrho^4=\varrho^2+\varrho$ and $\varrho^5=\varrho^4+1$.

\begin{proposition}\label{prop:plastic}
\begin{enumerate}[label=\textup{(\alph*)}]
\item The set $R_4=\{0,\varrho^2,\varrho^4,\varrho^5\}$ is a Sidon set whose six differences are exactly
  $1,\varrho,\varrho^2,\varrho^3,\varrho^4,\varrho^5$, each once.
\item The set $P_5=\{0,\varrho,\varrho^2,\varrho^3,\varrho^4\}$ has all ten positive differences in
  $\{\varrho^n:n\in\Z\}$; they are $\varrho^k$ for $k=-3,-2,-1,0,1,1,2,2,3,4$.
\item The set $R_5=\{0,\varrho^2,\varrho^4,\varrho^5,\varrho^7\}=\varrho^7-\varrho^3P_5$ has all ten positive differences in
  $\{\varrho^n:n\in\Z\}$; they are $\varrho^k$ for $k=0,1,2,3,4,4,5,5,6,7$.
\end{enumerate}
\end{proposition}

\begin{proof}
From $\varrho^3=\varrho+1$: $\varrho^4=\varrho^2+\varrho$, $\varrho^5=\varrho^3+\varrho^2=\varrho^2+\varrho+1=\varrho^4+1$,
$\varrho^6=\varrho^4+\varrho^3$ and $\varrho^7=\varrho^5+\varrho^4$. The differences of $R_5$ are
$\varrho^2,\varrho^4,\varrho^5,\varrho^7$ (from $0$); $\varrho^4-\varrho^2=\varrho$; $\varrho^5-\varrho^2=\varrho+1=\varrho^3$;
$\varrho^7-\varrho^2=\varrho^5+\varrho^4-\varrho^2=\varrho^2+2\varrho+1=\varrho^4+\varrho^3=\varrho^6$; $\varrho^5-\varrho^4=1$;
$\varrho^7-\varrho^4=\varrho^5$; $\varrho^7-\varrho^5=\varrho^4$. Those not involving $\varrho^7$ are the differences
of $R_4$: $\varrho^2,\varrho^4,\varrho^5,\varrho,\varrho^3,1$, which are distinct since $\varrho>1$.
Finally $\varrho^7-\varrho^3\cdot\{0,\varrho,\varrho^2,\varrho^3,\varrho^4\}=\{\varrho^7,\varrho^7-\varrho^4,\varrho^7-\varrho^5,
\varrho^7-\varrho^6,0\}=\{\varrho^7,\varrho^5,\varrho^4,\varrho^2,0\}$ by the relations above (note
$\varrho^7-\varrho^6=\varrho^5+\varrho^4-\varrho^4-\varrho^3=\varrho^5-\varrho^3=\varrho^2$), so $P_5$ is the image of $R_5$ under
$y\mapsto\varrho^{-3}(\varrho^7-y)$ and its differences are those of $R_5$ multiplied by $\varrho^{-3}$.
\end{proof}

\begin{remark}\label{rem:P5Lenstra}
The set $P_5$ is the analogue of $\{0,\varphi,\varphi^2,\varphi^3\}$. In \cite[(3.10)]{Lenstra1977} the bound $M\ge5$ is attributed to
\cite[(2.4)(c)]{Lenstra1977}: the sequence $0,1,x,1/(1-x),(x-1)/x$ has unit differences whenever $f(0),f(1),f(\zeta_6)$ are
units, $f$ the minimal polynomial of $x$ and $\zeta_6$ a primitive sixth root of unity. Among $\pm\varrho^{\pm1}$ this holds
exactly for $x=-\varrho$ ($f=X^3-X+1$, $f(\zeta_6)=-\zeta_6$) and $x=-\varrho^{-1}$ ($f=X^3-X^2+1$, $f(\zeta_6)=-\zeta_6^2$); for
$x=\varrho^{\pm1}$, $f(\zeta_6)\in\{-2-\zeta_6,\ \zeta_6^2-2\}$ has norm $7$. By $\varrho^3=\varrho+1$ the two sequences are
\[
  \{-\varrho,0,\varrho^{-3},1,\varrho^2\}=-\varrho+P_5,\qquad \{-\varrho^{-1},0,\varrho^{-2},1,\varrho^3\}=-\varrho^{-1}+\varrho^{-3}R_5 .
\]
Thus, up to maps $y\mapsto\pm\varrho^ky+t$, $P_5$ and $R_5$ are the sequences of \cite[(2.4)(c)]{Lenstra1977} for these
generators. As printed, \cite[(3.10)]{Lenstra1977} applies \cite[(2.4)(c)]{Lenstra1977} to $x=\alpha$, $\alpha^3-\alpha-1=0$, i.e.\ $x=\varrho$; for this $x$ the hypothesis of (2.4)(c) fails, since $f(\zeta_6)=-2-\zeta_6$ has norm~$7$. The printed $x=\alpha$ does satisfy the hypothesis of \cite[(2.4)(d)]{Lenstra1977}: $f(0)=f(1)=f(-1)=-1$ and $f(\varphi)=\varphi$, $f(\bar\varphi)=\bar\varphi$ (Lenstra's $\theta$ is a zero of $X^2-X-1$); as $\alpha+1=\alpha^3$, the sequence $0,1,\alpha,\alpha+1,\alpha^2$ of (2.4)(d) is $\varrho^{-1}P_5$, just as $\{0,\varphi,\varphi^2,\varphi^3\}$ is, up to dilation, the sequence of \cite[(2.4)(b)]{Lenstra1977} for $x=\varphi$. Alternatively, $x=-\alpha$ satisfies the hypothesis of (2.4)(c). Either way $M\ge5$, so the value $M=5$ in \cite[(3.10)]{Lenstra1977} is correct; only the printed reference to the part of (2.4) is affected. Independently, Nagell \cite[Th\'eor\`eme~2]{Nagell1964} lists all relations $1+E+E_1=0$ in units of the cubic field generated by a zero of $x^3+x^2-1$, i.e.\ by $\varrho^{-1}$. The hypothesis of \cite[(2.4)(c)]{Lenstra1977} and the
resulting set are invariant under $x\mapsto1/(1-x)$, whose orbits here are $\{-\varrho,\varrho^{-3},\varrho^{2}\}$ and
$\{-\varrho^{-1},\varrho^{-2},\varrho^{3}\}$; the second is the image of the first under $x\mapsto1-x$, which maps the set
to its image under $y\mapsto1-y$. Hence each of the six generators in these two orbits gives, up to maps $y\mapsto\pm\varrho^k y+t$, the
set $P_5$.
\end{remark}

Thus the analogue of Theorem~2$'$ fails for $\varrho$: five points suffice, and $\Dl(R_5)$ is even the
eight-term geometric progression $\{1,\varrho,\dots,\varrho^7\}$. By Remark~\ref{rem:lenstra}, $N(\varrho)=5$; thus $P_5$ and $R_5$ are extremal. For contrast, for $x=2$ at most three points are
possible: $2^p+2^q$ is a power of $2$ only if $p=q$, so four points would have three equal gaps $g$ and span
$3g$, which is not a power of $2$; and $\{0,1,2\}$ works. If $x$ is transcendental, already three points are
impossible, since $x^p+x^q=x^r$ is a non-trivial polynomial relation.

\subsection{What is classical and what is new}
\begin{itemize}
\item \emph{Classical or folklore:} the set $\{0,\varphi,\varphi^2,\varphi^3\}$ itself, i.e.\ the endpoints and the
  crossing points on a side of the regular pentagram (Euclid, \emph{Elements}~XIII.8; see Remark~\ref{rem:euclid}), which up to dilation is also
  the sequence $0,1,\varphi,\varphi^2$ of \cite[(2.4)(b), (3.9)]{Lenstra1977}, together with its difference multiset (Theorem~\ref{thm:T1}, a one-line consequence of
  $\varphi^2=\varphi+1$); Lemma~\ref{lem:phi}, Lemma~\ref{lem:cube}, Cassini's identity, the almost periodicity of
  trigonometric polynomials (Lemma~\ref{lem:trigAP}) and the periodicity criterion
  (Proposition~\ref{prop:periodic}), of which Theorem~\ref{thm:T4} is a direct consequence.
\item \emph{Known; we include a self-contained elementary proof:} the bound of four points, Theorem~\ref{thm:T2}(i) and
  Theorem~2$'$(iii), i.e.\ $N(\varphi)=4$; it follows from Lenstra \cite[(2.2), (3.9)]{Lenstra1977}, where the
  maximal sequence $0,1,\varphi,\varphi^2$ (a dilate of the golden quadruple) appears as the case $x=\varphi$ of
  \cite[(2.4)(b)]{Lenstra1977}, see \cite[(3.9)]{Lenstra1977}. Likewise $N(\varrho)=5$ follows from \cite[(2.2), (3.10)]{Lenstra1977}
  (Remark~\ref{rem:lenstra}). Lemma~\ref{lem:unit} and Theorem~2$'$(i) follow from Nagell's determination of
  the exceptional units \cite{Nagell1964,Nagell1969a,Nagell1969b} (Remark~\ref{rem:nagell}). Version~1 wrongly listed the
  four-point bound as new.
\item \emph{Explicit geometric formulation, own elementary proof:} the uniqueness results Theorem~\ref{thm:T2}(ii) and
  Theorem~2$'$(ii). They also follow from Nagell's list \cite{Nagell1964,Nagell1969a,Nagell1969b} after normalising an exceptional
  sequence (cf.\ \cite{LeutbecherMartinet1982}), see Remark~\ref{rem:nagell}; we do not claim them as new.
\item \emph{New (elementary):} the rigidity Theorem~3$'$, an explicit geometric rigidity statement with a self-contained
  elementary proof (it can also be derived from \cite[Lemma~3]{CampolinaMailybaev2021} by a short argument,
  Remark~\ref{rem:CM}; we have not found it in the literature); the sharp asymptotics of Proposition~\ref{prop:sharp} (the reduction of the shift norm to the
  maximum computed in \cite{Goldobin2007}, the exact value of $E_n$, its expansion and explicit remainder); and the order of
  $L(\varepsilon)$ in Proposition~\ref{prop:Leps}, with explicit constants for $0<\varepsilon\le\frac12$. The closed form
  (Proposition~\ref{prop:closed}) and the almost periods $\Fib_n/u$ (Proposition~\ref{prop:fibAP}) are elementary and
  serve as their basis; we have not found these statements for $C_S$ in the literature, and our own contribution there is
  the explicit value of $C_S(\Fib_n/u)$.
\item \emph{Supplementary:} the closure residues with their Fibonacci values (Theorem~\ref{thm:T3},
  Proposition~\ref{prop:fibres}), which concern the special quadruples $S_x$ and are elementary companions of
  Theorem~3$'$, and the explicit plastic-number sets (Proposition~\ref{prop:plastic}).
\item \emph{Not claimed:} depth beyond the elementary level, or any meaning of the scale $K$ other than a
  normalisation.
\end{itemize}

\subsection{Open questions}
\begin{question}\label{q:N}
Let $N(x)$ be as in Definition~\ref{def:N}. We know $N(\varphi)=4$ and $N(\varrho)=5$ (Lenstra \cite[(3.9), (3.10)]{Lenstra1977}; Theorem~2$'$,
Proposition~\ref{prop:plastic}, Remark~\ref{rem:lenstra}) and $N(2)=3$. Determine $N(x)$ for every $x>1$. Is
$\sup_{x>1}N(x)=5$? Since $N(x)\ge3$ if and only if $x^l=x^k+1$ for some integers $l>k\ge0$, a key step is to decide
when two distinct trinomials $X^l-X^k-1$ have a common root $x>1$. For $x\ge1.05$ this was settled, with a
computer-assisted proof, by Campolina and Mailybaev \cite[Lemma~3]{CampolinaMailybaev2021}, who conjecture that their classification of nondegenerate logarithmic lattices (their Theorem~1) holds for all $x>1$ ($\lambda>1$ in their notation), with no
bound on the exponents \cite[Remark~2]{CampolinaMailybaev2021}; the finite list $\varrho^{1/m}$, $m\le5$, of
\cite[Lemma~3]{CampolinaMailybaev2021} belongs to the range $x\ge1.05$. A relevant structural tool is
the irreducibility theorem of Ljunggren \cite{Ljunggren1960} and Tverberg \cite{Tverberg1960}, preceded by
Selmer's work on $X^n-X-1$ \cite{Selmer1956}: apart from cyclotomic factors, $X^n\pm X^m\pm1$ is irreducible; see
\cite[Lemma~4.1]{BoydMartinThom2015} for a statement covering exponents with $\gcd(n,m)>1$. Aarts, Fokkink and
Kruijtzer \cite{AartsFokkinkKruijtzer2001} classify the common positive roots of $X^n-X-1$ and $X^m-X^{m-1}-1$; for two
distinct polynomials the only one is $\varrho$. Since $X^l-X^k-1$ has a root $x>1$, which is not a root of unity, its
non-cyclotomic part is the minimal polynomial of $x$ (for $k=0$, $X^l-2$ is irreducible by Eisenstein's criterion). Hence two such trinomials have a common root $x>1$ if and
only if their non-cyclotomic parts coincide. The author will treat Question~\ref{q:N} in a separate paper.
\end{question}

\begin{question}
Classify the triples $(x,m,k)$ for which there is an $m$-point set $A\subset\R$ whose set of distinct differences is
exactly a geometric progression $\{c,cx,\dots,cx^k\}$. Theorem~3$'$ settles $(m,k)=(4,3)$ (only $x=\varphi$), while
Proposition~\ref{prop:plastic} gives examples with $(m,k)=(4,5)$ and $(5,7)$ for $x=\varrho$.
\end{question}

\begin{question}
By Proposition~\ref{prop:Leps}, $L(\varepsilon)\asymp1/\varepsilon$. Does $\lim_{\varepsilon\to0}\varepsilon L(\varepsilon)$ exist,
and what is its value (for $u=1$)? The constants in Proposition~\ref{prop:Leps} differ by a factor of about $443$ ($6\pi\varphi^3\approx79.8$ versus $\sqrt2/(3\varphi^2)\approx0.18$).
\end{question}

\section{Formal verification}\label{sec:formal}
\setlength{\emergencystretch}{3em}%

The combinatorial and algebraic core of Sections~\ref{sec:T1}--\ref{sec:T4} is formalised in the proof
assistant Lean~4 \cite{deMouraUllrich2021} on top of the Mathlib library \cite{Mathlib2020};
Table~\ref{tab:lean} lists exactly which results are covered. The Bohr almost periodicity of $C_S$
(Lemma~\ref{lem:trigAP}, Proposition~\ref{prop:periodic}, Theorem~\ref{thm:T4}(b)) and
Proposition~\ref{prop:plastic} are not formalised. The formalisation is an independent check of part of
the proofs given above; the paper is self-contained and does not depend on it.

\subsection{Environment and reproduction}
\begin{itemize}
\item Lean version \texttt{v4.35.0}; the file \texttt{lean-toolchain} contains
  \texttt{leanprover/lean4:v4.35.0}.
\item Mathlib: tag \texttt{v4.35.0}, commit \texttt{56a005855ea17e0b0958f73ec801ae140b0a1679}, pinned in
  \texttt{lake-manifest.json}.
\item Lake project \texttt{GoldenDifference}; all declarations lie in the namespace
  \texttt{GoldenDifference}. Modules: \texttt{Basic}, \texttt{Extremality}, \texttt{Rigidity},
  \texttt{Algebra}, \texttt{Orbit}, \texttt{Fibonacci}, \texttt{Cosine}, \texttt{AlmostPeriods}; audit file
  \texttt{Audit.lean}.
\item Source code: to be deposited in a public archive, with a persistent identifier and the commit hash,
  so that Table~\ref{tab:lean} can be checked independently; until then it is available from the author on
  request.
\end{itemize}
From the root directory of the project, the formalisation is checked with the commands
\begin{verbatim}
lake exe cache get        # fetch precompiled Mathlib files
lake build                # compile and check all modules
lake env lean Audit.lean  # print the axioms of the audited declarations
\end{verbatim}

\subsection{Conventions}
\begin{itemize}
\item A finite set $A\subset\R$ is a term of type \texttt{Finset}~$\R$. Its positive difference set
  $\Dl(A)$ (Definition~\ref{def:diff}) is \texttt{posDiff A}, the set of all $x-y$ with $x,y\in A$ and
  $y<x$. The difference multiset is \texttt{differenceMultiset A}, which has one entry $x-y$ for each pair
  $y<x$ in $A$; the multiplicity of $d$ in it is $r_A(d)$.
\item $\varphi$ is Mathlib's \texttt{goldenRatio}, defined as $(1+\sqrt5)/2$. The set \texttt{orbit c} is
  $\Phi_c=\{c\varphi^n:n\in\Z\}$ of Theorem~2$'$, with integer powers.
\item \texttt{cosineSum u t} is defined as the right-hand side
  $\cos(2\pi ut)+2\cos(2\pi\varphi ut)+2\cos(2\pi\varphi^2ut)+\cos(2\pi\varphi^3ut)$; its identification with
  $C_S(t)$ of Definition~\ref{def:cosine}, which follows from Theorem~\ref{thm:T1}, is proved on paper and is not
  formalised. By Remark~\ref{rem:structure}(d) this expression depends only on $u$.
\item The Fibonacci numbers $\Fib_n$ are Mathlib's \texttt{Nat.fib n}, with $\Fib_0=0$ and $\Fib_1=1$ as in
  Section~\ref{sec:prelim}.
\item $\nint{\theta}$ is \texttt{distanceToInteger}~$\theta=|\theta-\operatorname{round}(\theta)|$, which equals
  $\min_{m\in\Z}|\theta-m|$.
\item Bohr almost periodicity (Definition~\ref{def:ap}) is not defined in the project.
\end{itemize}

\subsection{Formalised statements}
Table~\ref{tab:lean} lists the results of the paper, the corresponding declarations in the namespace
\texttt{GoldenDifference}, and the status of each result: \emph{formalised} (the listed declarations prove
the result, with the qualifications stated in the table), \emph{partly formalised} (only the parts named in
the table are proved in Lean), or \emph{not formalised}.

\begingroup
\small
\begin{longtable}{@{}
  >{\raggedright\arraybackslash}p{0.27\textwidth}
  >{\raggedright\arraybackslash}p{0.33\textwidth}
  >{\raggedright\arraybackslash}p{0.34\textwidth}@{}}
\caption{Results of the paper, the corresponding Lean declarations (namespace \texttt{GoldenDifference}),
and the status of their formalisation.}\label{tab:lean}\\
\toprule
result & Lean declaration(s) & status\\
\midrule
\endfirsthead
\toprule
result & Lean declaration(s) & status\\
\midrule
\endhead
\bottomrule
\endfoot
Theorem~\ref{thm:T1} (six differences, $\Dl(S)=\Dset$, multiplicities $1,2,2,1$)
  & \texttt{six\_identities}, \texttt{theorem1\_set}, \texttt{theorem1\_multiset},
    \texttt{theorem1\_multiplicities}
  & formalised, for every $F\in\R$ and $u>0$\\
\addlinespace
Lemma~\ref{lem:cube}, Remark~\ref{rem:structure}(a)
  & \texttt{posDiff\_four}, \texttt{differenceMultiset\_four}, \texttt{canonical\_eq\_reflected},
    \texttt{canonical\_card}
  & partly formalised: the differences of four increasing reals; the translate form of $S$\\
\addlinespace
Theorem~\ref{thm:T2}(i), (ii)
  & \texttt{theorem2\_card}, \texttt{theorem2\_unique}
  & formalised\\
\addlinespace
Lemma~\ref{lem:unit} (golden unit equation)
  & \texttt{golden\_unit}, \texttt{golden\_unit\_set}
  & formalised\\
\addlinespace
Theorem~2$'$
  & \texttt{orbit\_card}
  & partly formalised: only (iii), at most four points with all differences in $\Phi_c$; (i), (ii) not
    formalised\\
\addlinespace
Theorem~\ref{thm:T3}
  & \texttt{base\_unique}, \texttt{closure\_residues}, \texttt{closure\_iff}
  & partly formalised: the residue identities and $\rho_-(x)=0\iff x=\varphi$; (a), (c), (d) and the
    equivalence with $\Dl(S_x)=\Dset_x$ in (b) not formalised\\
\addlinespace
Lemma~\ref{lem:cassini}, Proposition~\ref{prop:fibres}
  & \texttt{cassini}, \texttt{fibonacci\_residue}, \texttt{fibonacci\_closure\_residues}
  & formalised: Cassini's identity in the equivalent form
    $\Fib_{n+1}^2-\Fib_{n+1}\Fib_n-\Fib_n^2=(-1)^n$; the residue values, for $n\ge1$ (the identities are
    algebraic; Proposition~\ref{prop:fibres} takes $n\ge2$ only because $x_1=1$ is not admissible in
    Definition~\ref{def:Sx})\\
\addlinespace
Theorem~3$'$ (rigidity)
  & \texttt{theorem3\_prime}, \texttt{theorem3\_prime\_iff}
  & formalised: the equivalence, and $T$ is a translate of $c\{0,\varphi,\varphi^2,\varphi^3\}$\\
\addlinespace
Lemma~\ref{lem:phi}(iv) (irrationality, incommensurability)
  & \texttt{phi\_irrational}, \texttt{no\_common\_integer\_multiple}
  & formalised\\
\addlinespace
Proposition~\ref{prop:closed} (closed form, (i)--(iii))
  & \texttt{cosine\_closed\_form}, \texttt{cosine\_le\_six}, \texttt{cosine\_lower\_bound},
    \texttt{cosine\_eq\_six\_iff}, \texttt{cosine\_minimum\_iff}, \texttt{cosine\_minimizers}
  & formalised: the closed form $8\cos^2(\pi at)\cos^2(\pi bt)-2$ for \texttt{cosineSum} and (i)--(iii); the
    identity for $|\widehat S(t)|^2$ not formalised\\
\addlinespace
Theorem~\ref{thm:T4}(a) (non-periodicity)
  & \texttt{cosine\_not\_periodic}
  & formalised\\
\addlinespace
Proposition~\ref{prop:fibAP}(i), (ii)
  & \texttt{fibonacci\_error}, \texttt{fibonacci\_error\_abs}, \texttt{fibonacci\_nearest\_distance},
    \texttt{cosine\_shift\_bound}, \texttt{fibonacci\_almost\_period}
  & formalised (the equality in (i) for $n\ge2$); (iii) not formalised\\
\addlinespace
Lemma~\ref{lem:trigAP}, Proposition~\ref{prop:periodic}, Theorem~\ref{thm:T4}(b) (Bohr almost periodicity)
  & ---
  & not formalised\\
\addlinespace
Proposition~\ref{prop:plastic} (plastic number)
  & ---
  & not formalised\\
\end{longtable}
\endgroup

\subsection{What is not formalised}
The following results are proved by hand in the paper and are not formalised:
\begin{itemize}
\item Theorem~2$'$(i), (ii);
\item parts (a), (c) and (d) of Theorem~\ref{thm:T3}, and in part (b) the equivalence with
  $\Dl(S_x)=\Dset_x$;
\item Lemma~\ref{lem:trigAP}, Proposition~\ref{prop:periodic} and Theorem~\ref{thm:T4}(b), i.e.\ the Bohr
  almost periodicity of $C_S$; note that the formalised non-periodicity (Theorem~\ref{thm:T4}(a)) does not by
  itself give almost periodicity;
\item Lemma~\ref{lem:square} and the identity $|\widehat S(t)|^2=4+2C_S(t)=16\cos^2(\pi at)\cos^2(\pi bt)$
  in Proposition~\ref{prop:closed};
\item Proposition~\ref{prop:fibAP}(iii);
\item Proposition~\ref{prop:plastic};
\item Lemma~\ref{lem:dense}, Proposition~\ref{prop:sharp} and Proposition~\ref{prop:Leps} (new in version~2);
\item Corollary~\ref{cor:ladder} and Lemma~\ref{lem:ladder} (the golden ladder);
\item Remark~\ref{rem:torus} (torus representation).
\end{itemize}
Neither are the numerical values of Example~\ref{ex:planck} and Section~\ref{sec:numerics}, which are
floating-point computations. More generally, no result of the paper is claimed to be formally
checked unless it appears in Table~\ref{tab:lean} as formalised or partly formalised, and in the latter case
only to the extent stated there.

\subsection{Results}
In a build carried out on 30 October 2026 in the environment above, \texttt{lake build} terminated without
errors, with the message \texttt{Build completed successfully (2433 jobs)}. The sources of the project
contain no occurrence of \texttt{sorry}, \texttt{admit}, \texttt{axiom}, \texttt{native\_decide} or
\texttt{opaque}. Hence every declaration of the project, in particular every declaration listed in
Table~\ref{tab:lean}, has been accepted by the Lean kernel with a complete proof, and the project adds no
axioms of its own.

The file \texttt{Audit.lean} lists the $27$ audited declarations and runs \texttt{\#print axioms} on each of
them. For every one of these $27$ declarations the output contains no axioms other than the three standard
axioms of Lean, \texttt{propext}, \texttt{Classical.choice} and \texttt{Quot.sound}; in particular
\texttt{sorryAx} does not occur. The remaining eight declarations of Table~\ref{tab:lean} are not audited
directly; each of them is used in the proof of an audited declaration.

\clearpage
\section*{Changes in version 2}
\addcontentsline{toc}{section}{Changes in version 2}
\begin{itemize}
\item Priority corrected: the bound of four points in Theorems~\ref{thm:T2}(i) and~2$'$(iii) ($N(\varphi)=4$) and
  $N(\varrho)\le5$ follow from Lenstra~\cite{Lenstra1977}, where the sequence $0,1,\varphi,\varphi^2$ also appears
  (Remark~\ref{rem:lenstra}); Lemma~\ref{lem:unit} and Theorem~2$'$(i) follow from Nagell's exceptional units
  (Remark~\ref{rem:nagell}); Theorem~\ref{thm:T4} is listed as a consequence of classical facts. Uniqueness
  (Theorems~\ref{thm:T2}(ii), 2$'$(ii)) is presented as an explicit geometric formulation with an elementary proof (it
  also follows from Nagell via normalisation); rigidity (Theorem~3$'$) remains new as an explicit statement with a self-contained proof; it can also be derived from \cite[Lemma~3]{CampolinaMailybaev2021} by a short argument (Remark~\ref{rem:CM}).
\item The computational observation on six-point sets for $\varrho$ is replaced by the theorem $N(\varrho)=5$;
  Proposition~\ref{prop:plastic} now includes $P_5=\{0,\varrho,\dots,\varrho^4\}$.
\item References to Nagell \cite{Nagell1964}, Naito \cite{Naito1996}, Aarts--Fokkink--Kruijtzer,
  Boyd--Martin--Thom, Selmer, Tverberg \cite{AartsFokkinkKruijtzer2001,BoydMartinThom2015,Selmer1956,Tverberg1960}, Moffet
  \cite{Moffet1968} and Wei\ss{} \cite{Weiss2002}; the short derivation of Theorem~3$'$ from
  \cite[Lemma~3]{CampolinaMailybaev2021} (Remark~\ref{rem:CM}).
\item References to Goldobin--Koelle--Kleiner--Buzdin \cite{Goldobin2007} (Proposition~\ref{prop:sharp}), Anikushin, Goda and
  Slater \cite{Anikushin2017,Goda2024,Slater1967} (Proposition~\ref{prop:Leps}) and Leutbecher--Martinet
  \cite{LeutbecherMartinet1982}.
\item Definition~\ref{def:N} of $N(x)$; Question~\ref{q:N} updated; references to Campolina and
  Mailybaev~\cite{CampolinaMailybaev2021}, Ljunggren~\cite{Ljunggren1960} and Nagell~\cite{Nagell1969a,Nagell1969b}.
\item Added: Lemma~\ref{lem:dense} (density, replacing the appeal to Kronecker's theorem),
  Proposition~\ref{prop:sharp} (sharp constants for the almost periods $\Fib_n/u$) and
  Proposition~\ref{prop:Leps} ($L(\varepsilon)\asymp1/\varepsilon$, answering a question of version~1).
\item Minor corrections; abstract and introduction revised accordingly.
\end{itemize}

\section*{Changes in version 3}
\addcontentsline{toc}{section}{Changes in version 3}
\begin{itemize}
\item Remark~\ref{rem:lenstra}: explicit reference to Leutbecher--Martinet \cite[\S3.0, Thm.~4.1.1]{LeutbecherMartinet1982} for $N(\varphi)=4$ and $N(\varrho)=5$.
\item DOIs added to the bibliography; sampling remark in Section~\ref{sec:numerics} clarified.
\end{itemize}

\section*{Acknowledgements}
\addcontentsline{toc}{section}{Acknowledgements}
The author used AI tools (Claude, Anthropic; Codex, OpenAI) for drafting and editing the text and for the
Lean formalisation. All mathematical content was verified by the author. Theorems~\ref{thm:T1},
\ref{thm:T2} and~3$'$, Lemmas~\ref{lem:phi}(iv), \ref{lem:unit} and~\ref{lem:cassini},
Theorem~2$'$(iii), Proposition~\ref{prop:fibres}, Proposition~\ref{prop:closed} (except the identity for
$|\widehat S|^2$), Theorem~\ref{thm:T4}(a) and
Proposition~\ref{prop:fibAP}(i)--(ii) were in addition formally checked in Lean~4, as were parts of
Lemma~\ref{lem:cube}, Remark~\ref{rem:structure}(a) and Theorem~\ref{thm:T3}, in each case to the extent
stated in Table~\ref{tab:lean} (Section~\ref{sec:formal}); the remaining results were verified by hand.

\bibliographystyle{plainnat}
\bibliography{refs}

\end{document}